\documentclass[a4paper,12pt]{article}

\usepackage{amsfonts,amssymb,amsmath,amsthm,latexsym,epsf,epsfig,amscd,bbm,stmaryrd}

\usepackage[english]{babel}
\makeatletter \@addtoreset{equation}{section} 
\makeatother

\makeatletter \@addtoreset{enunciato}{section} 
\makeatother

\newcounter{enunciato}[section]

\newtheorem{ittheorem}{Theorem}
\newtheorem{itlemma}{Lemma}
\newtheorem{itproposition}{Proposition}
\newtheorem{itdefinition}{Definition}
\newtheorem{itremark}{Remark}
\newtheorem{itclaim}{Claim}
\newtheorem{itfact}{Fact}
\newtheorem{itconjecture}{Conjecture}
\newtheorem{itcorollary}{Corollary}

\newenvironment{theorem}{\addtocounter{enunciato}{1}
\begin{ittheorem}}{\end{ittheorem}}

\newenvironment{lemma}{\addtocounter{enunciato}{1}
\begin{itlemma}}{\end{itlemma}}

\newenvironment{proposition}{\addtocounter{enunciato}{1}
\begin{itproposition}}{\end{itproposition}}

\newenvironment{definition}{\addtocounter{enunciato}{1}
\begin{itdefinition}}{\end{itdefinition}}

\newenvironment{remark}{\addtocounter{enunciato}{1}
\begin{itremark}}{\end{itremark}}

\newenvironment{conjecture}{\addtocounter{enunciato}{1}
\begin{itconjecture}}{\end{itconjecture}}

\newenvironment{corollary}{\addtocounter{enunciato}{1}
\begin{itcorollary}}{\end{itcorollary}}

\newcommand{\be}[1]{\begin{equation}\label{#1}}
\newcommand{\ee}{\end{equation}}

\newcommand{\bl}[1]{\begin{lemma}\label{#1}}
\newcommand{\el}{\end{lemma}}

\newcommand{\br}[1]{\begin{remark}\label{#1}}
\newcommand{\er}{\end{remark}}

\newcommand{\bt}[1]{\begin{theorem}\label{#1}}
\newcommand{\et}{\end{theorem}}

\newcommand{\bd}[1]{\begin{definition}\label{#1}}
\newcommand{\ed}{\end{definition}}

\newcommand{\bp}[1]{\begin{proposition}\label{#1}}
\newcommand{\ep}{\end{proposition}}

\newcommand{\bc}[1]{\begin{corollary}\label{#1}}
\newcommand{\ec}{\end{corollary}}

\newcommand{\bcj}[1]{\begin{conjecture}\label{#1}}
\newcommand{\ecj}{\end{conjecture}}

\def \Z {{\mathbb Z}}
\def \R {{\mathbb R}}

\def \N {{\mathbb N}}

\def \ba {\begin{array}}
\def \ea {\end{array}}

\def \P  {{\mathbb P}}
\def \E  {{\mathbb E}}
\def \H {{\mathbb H}}
\def \cP {{\mathcal P}}

\def \cH {{\mathcal H}}
\def \cG {{\mathcal G}}
\def \cF {{\mathcal F}}

\def \erom  {\mathrm{e}}
\def \di {\mathrm{d}}
\def \trn {{\interleave}}

\begin{document}


\title{Law of large numbers for a class of\\ 
random walks in dynamic random environments}

\author{\renewcommand{\thefootnote}{\arabic{footnote}}
L.\ Avena \footnotemark[1]
\\
\renewcommand{\thefootnote}{\arabic{footnote}}
F.\ den Hollander \footnotemark[1]\,\,\,\,\footnotemark[2]
\\
\renewcommand{\thefootnote}{\arabic{footnote}}
F.\ Redig \footnotemark[1]}

\footnotetext[1]{
Mathematical Institute, Leiden University, P.O.\ Box 9512, 
2300 RA Leiden, The Netherlands} 

\footnotetext[2]{
EURANDOM, P.O.\ Box 513, 5600 MB Eindhoven, The Netherlands}

\maketitle

\begin{abstract}
In this paper we consider a class of one-dimensional interacting particle systems in 
equilibrium, constituting a dynamic random environment, together with a nearest-neighbor 
random walk that on occupied/vacant sites has a local drift to the right/left. We adapt 
a regeneration-time argument originally developed by Comets and Zeitouni~\cite{CoZe04} 
for static random environments to prove that, under a space-time mixing property for 
the dynamic random environment called cone-mixing, the random walk has an a.s.\ constant 
global speed. In addition, we show that if the dynamic random environment is exponentially 
mixing in space-time and the local drifts are small, then the global speed can be written 
as a power series in the size of the local drifts. From the first term in this series the 
sign of the global speed can be read off.

The results can be easily extended to higher dimensions. 

\vspace{0.5cm}\noindent
{\it Acknowledgment.} The authors are grateful to R.\ dos Santos and V.\ Sidoravicius
for fruitful discussions.

\end{abstract}

\newpage


\section{Introduction and main result}
\label{S1} 

In Section~\ref{S1} we define the random walk in dynamic random environment, introduce a 
space-time mixing property for the random environment called cone-mixing, and state our 
law of large numbers for the random walk subject to cone-mixing. In Section~\ref{S2} we 
give the proof of the law of large numbers with the help of a space-time regeneration-time
argument. In Section~\ref{S3} we assume a stronger space-time mixing property, namely, 
exponential mixing, and derive a series expansion for the global speed of the random walk 
in powers of the size of the local drifts. This series expansion converges for small enough 
local drifts and its first term allows us to determine the sign of the global speed. (The 
perturbation argument underlying the series expansion provides an alternative proof of the 
law of large numbers.) In Appendix~\ref{appA} we give examples of random environments that 
are cone-mixing. In Appendix~\ref{appB} we compute the first three terms in the expansion 
for an independent spin-flip dynamics.


\subsection{Model}
\label{S1.1}

Let $\Omega = \{0,1\}^\Z$. Let $C(\Omega)$ be the set of continuous functions on 
$\Omega$ taking values in $\R$, $\cP(\Omega)$ the set of probability measures on 
$\Omega$, and $D_\Omega[0,\infty)$ the path space, i.e., the set of c\`adl\`ag 
functions on $[0,\infty)$ taking values in $\Omega$. In what follows, 
\be{xidef}
\xi = (\xi_t)_{t \geq 0} \quad \mbox{ with } \quad
\xi_t = \{\xi_t(x)\colon\,x\in\Z\}
\ee
is an interacting particle system taking values in $\Omega$, with $\xi_t(x)=0$ meaning 
that site $x$ is vacant at time $t$ and $\xi_t(x)=1$ that it is occupied. The paths of 
$\xi$ take values in $D_\Omega[0,\infty)$. The law of $\xi$ starting from $\xi_0=\eta$ 
is denoted by $P^\eta$. The law of $\xi$ when $\xi_0$ is drawn from $\mu\in\cP(\Omega)$ 
is denoted by $P^\mu$, and is given by 
\be{PathMeas}
P^\mu(\cdot) = \int_{\Omega} P^\eta(\cdot)\,\mu(\di\eta).
\ee
Through the sequel we will assume that 
\be{Pmuass}
P^\mu \mbox{ is stationary and ergodic under space-time shifts}.
\ee 
Thus, in particular, $\mu$ is a homogeneous extremal equilibrium for $\xi$. The Markov 
semigroup associated with $\xi$ is denoted by $S_\mathrm{IPS}=(S_\mathrm{IPS}(t))_{t \geq 0}$. 
This semigroup acts from the left on $C(\Omega)$ as  
\be{Semigroup} 
\big(S_\mathrm{IPS}(t)f\big)(\cdot) = E^{(\cdot)}[f(\xi_t)], \qquad f\in C(\Omega),
\ee
and acts from the right on $\cP(\Omega)$ as
\be{dualSemigroup}
\big(\nu S_\mathrm{IPS}(t)\big)(\cdot) = P^\nu(\xi_t \in \cdot\,), \qquad \nu \in \cP(\Omega). 
\ee 
See Liggett~\cite{Li85}, Chapter I, for a formal construction.

Conditional on $\xi$, let
\be{rwdef}
X = \{X(t)\colon\, t\geq 0\}
\ee
be the random walk with local transition rates
\be{rwtrans}
\begin{aligned}
&x \to x+1 \quad \mbox{ at rate } \quad \alpha\,\xi_t(x) + \beta\,[1-\xi_t(x)],\\
&x \to x-1 \quad \mbox{ at rate } \quad \beta\,\xi_t(x) + \alpha\,[1-\xi_t(x)],
\end{aligned}
\ee
where w.l.o.g.
\be{albe}
0<\beta<\alpha<\infty.
\ee 
Thus, on occupied sites the random walk has a local drift to the right while on vacant sites 
it has a local drift to the left, of the same size. Note that the sum of the jump rates 
$\alpha+\beta$ is independent of $\xi$. Let $P^\xi_0$ denote the law of $X$ starting from 
$X_0=0$ conditional on $\xi$, which is the \emph{quenched} law of $X$. The \emph{annealed} 
law of $X$ is
\be{annealed} 
\P_{\mu,0}(\cdot)
= \int_{D_\Omega[0,\infty)} P^\xi_0(\cdot)\,P^{\mu}(\di\xi). 
\ee


\subsection{Cone-mixing and law of large numbers}
\label{S1.2}

In what follows we will need a \emph{mixing property} for the law $P^\mu$ of $\xi$.
Let $(\cdot,\cdot)$ and $\|\cdot\|$ denote the inner product, respectively, the 
Euclidean norm on $\R^2$. Put $\ell=(0,1)$. For $\theta \in (0,\tfrac12\pi)$ and 
$t \geq 0$, let
\be{Cone}
C_t^\theta 
= \big\{u\in\Z \times [0,\infty)\colon\,(u-t\ell,\ell) \geq \|u-t\ell\|\cos\theta\big\}
\ee 
be the cone whose tip is at $t\ell=(0,t)$ and whose wedge opens up in the direction 
$\ell$ with an angle $\theta$ on either side (see Figure~\ref{fig-cone}). Note that 
if $\theta=\tfrac12\pi$ ($\theta=\tfrac14\pi$), then the cone is the half-plane 
(quarter-plane) above $t\ell$.


\begin{figure}[hbtp]
\vspace{1.5cm}
\begin{center}
\setlength{\unitlength}{0.3cm}
\begin{picture}(12,10)(-8,0)
\put(-8,0){\line(16,0){16}} 
\put(0,0){\line(0,12){12}}
{\thicklines 
\qbezier(0,4)(3,7)(6,10) 
\qbezier(0,4)(-3,7)(-6,10)
}
\put(-1.3,-1.7){$(0,0)$}
\put(.8,3.7){$(0,t)$}
\put(-13,0){\vector(1,0){4}}
\put(-6,8){\vector(1,0){4}}
\put(1,0){\circle*{.5}}
\put(2,0){\circle*{.5}}
\put(3,0){\circle*{.5}}
\put(4,0){\circle*{.5}}
\put(5,0){\circle*{.5}}
\put(0,0){\circle*{.5}}
\put(-1,0){\circle*{.5}}
\put(-2,0){\circle*{.5}}
\put(-3,0){\circle*{.5}}
\put(-4,0){\circle*{.5}}
\put(-5,0){\circle*{.5}}
\put(-1.1,5.7){$\theta$}
\put(.5,5.7){$\theta$}
\put(7,6){$\Z \times [0,\infty)$}
\put(-14.5,-.3){$\Z$}
\put(-8.3,7.8){$C_t^\theta$}
\put(-1,12.8){\small\mbox{time}}
\put(9,-.3){\small\mbox{space}}
\end{picture}
\end{center}
\caption{\small The cone $C_t^\theta$.}
\label{fig-cone}
\vspace{0.3cm}
\end{figure}


\bd{PhiMixing}
A probability measure $P^{\mu}$ on $D_\Omega[0,\infty)$ satisfying {\rm (\ref{Pmuass})} 
is said to be cone-mixing if, for all $\theta \in (0,\tfrac12\pi)$,
\be{ConeMix1}
\lim_{t\to\infty} \sup_{ {A\in\cF_0,\,B\in\cF_t^\theta} \atop {P^\mu(A)>0} }
\Big|P^{\mu}(B\mid A)-P^{\mu}(B)\Big| = 0,
\ee 
where
\be{Fdefs}
\begin{aligned}
\cF_0 &= \sigma\big\{\xi_0(x) \colon\,x \in \Z\big\},\\
\cF_t^\theta &= \sigma\big\{\xi_s(x)\colon\,(x,s) \in C_t^\theta\big\}.
\end{aligned}
\ee
\ed 

\noindent
In Appendix~\ref{appA} we give examples of interacting particle systems that are 
cone-mixing. 

We are now ready to formulate our law of large numbers (LLN).

\bt{LLN} 
Assume {\rm (\ref{Pmuass})}. If $P^{\mu}$ is cone-mixing, then there exists a $v\in\R$ 
such that
\be{LLN1}
\lim_{t\to\infty} X_t/t = v \qquad \P_{\mu,0}-a.s.
\ee
\et

\noindent
The proof of Theorem~\ref{LLN} is given in Section~\ref{S2}, and is based on a 
\emph{regeneration-time argument} originally developed by Comets and Zeitouni 
\cite{CoZe04} for static random environments (based on earlier work by Sznitman 
and Zerner~\cite{SzZe99}).

We have no criterion for when $v<0$, $v=0$ or $v>0$. In view of (\ref{albe}), a 
naive guess would be that these regimes correspond to $\rho<\tfrac12$, $\rho=\tfrac12$ 
and $\rho>\tfrac12$, respectively, with $\rho=P^\mu(\xi_0(0)=1)$ the density of 
occupied sites. However, $v=(2\tilde\rho-1)(\alpha-\beta)$, with $\tilde\rho$ the 
asymptotic fraction of time spent by the walk on occupied sites, and the latter 
is a non-trivial function of $P^\mu$, $\alpha$ and $\beta$. We do not (!) expect that 
$\tilde\rho=\tfrac12$ when $\rho=\tfrac12$ in general. Clearly, if $P^\mu$ is invariant 
under swapping the states $0$ and $1$, then $v=0$.


\subsection{Global speed for small local drifts}
\label{S1.3}

For small $\alpha-\beta$, $X$ is a perturbation of simple random walk. In that case 
it is possible to derive an expansion of $v$ in powers of $\alpha-\beta$, provided 
$P^\mu$ satisfies an exponential space-time mixing property referred to as $M<\epsilon$ 
(Liggett~\cite{Li85}, Section I.3). Under this mixing property, $\mu$ is even uniquely 
ergodic.

Suppose that $\xi$ has shift-invariant local transition rates
\be{cTdef}
c(A,\eta), \qquad A\subset\Z \mbox{ finite},\,\eta\in\Omega,
\ee 
i.e., $c(A,\eta)$ is the rate in the configuration $\eta$ to change the states at the 
sites in $A$, and $c(A,\eta)=c(A+x,\tau_x\eta)$ for all $x\in\Z$ with $\tau_x$ the shift 
of space over $x$. Define
\be{Mepsdef}
\begin{aligned}
M &= \sum_{A \ni 0} \sum_{x \neq 0}\,\sup_{\eta\in\Omega}
|c(A,\eta)-c(A,\eta^x)|,\\
\epsilon &= \inf_{\eta\in\Omega} \sum_{A \ni 0} 
|c(A,\eta) + c(A,\eta^0)|,
\end{aligned}
\ee
where $\eta^x$ is the configuration obained from $x$ by changing the state at site $x$.
The interpretation of (\ref{Mepsdef}) is that $M$ is a measure for the \emph{maximal} 
dependence of the transition rates on the states of single sites, while $\epsilon$ is 
a measure for the \emph{minimal} rate at which the states of single sites change. See 
Liggett~\cite{Li85}, Section I.4, for examples.   
  
\bt{PERT}
Assume {\rm (\ref{Pmuass})} and suppose that $M<\epsilon$. If $\alpha-\beta<\tfrac12
(\epsilon-M)$, then 
\be{pert}
v = \sum_{n\in\N} c_n\,(\alpha-\beta)^n \in \R \quad 
\mbox{ with } \quad c_n=c_n(\alpha+\beta;P^\mu),
\ee
where $c_1=2\rho-1$ and $c_n\in\R$, $n \in\N\backslash\{1\}$, are given by a recursive 
formula (see Section {\rm \ref{S3.3}}).
\et

\noindent
The proof of Theorem~\ref{PERT} is given in Section~\ref{S3}, and is based on an analysis 
of the semigroup associated with the environment proces, i.e., the environment as seen 
relative to the random walk. The generator of this process turns out to be a sum of a 
large part and a small part, which allows for a perturbation argument. In Appendix~\ref{appA} 
we show that $M<\epsilon$ implies cone-mixing for spin-flip systems, i.e., systems for which
$c(A,\eta)=0$ when $|A| \geq 2$. 

It follows from Theorem~\ref{PERT} that for $\alpha-\beta$ small enough the global speed
$v$ changes sign at $\rho=\tfrac12$:
\be{weaklim}
v = (2\rho-1)(\alpha-\beta) + O\big((\alpha-\beta)^2\big) \mbox{ as } \alpha\downarrow\beta
\mbox{ for } \rho \mbox{ fixed}.
\ee
We will see in Section~\ref{S3.3} that $c_2=0$ when $\mu$ is a \emph{reversible} equilibrium,
in which case the error term in (\ref{weaklim}) is $O((\alpha-\beta)^3)$. 

In Appendix~\ref{appB} we consider an independent spin-flip dynamics such that $0$ changes 
to $1$ at rate $\gamma$ and $1$ changes to $0$ at rate $\delta$, where $0<\gamma,\delta<\infty$. 
By reversibility, $c_2=0$. We show that
\be{third}
c_3 = \frac{4}{U^2}\,\rho(1-\rho)(2\rho-1)\,f(U,V),
\quad f(U,V) = \frac{2U+V}{\sqrt{V^2+2UV}}-\frac{2U+2V}{\sqrt{V^2+UV}}+1,
\ee
with $U=\alpha+\beta$, $V=\gamma+\delta$ and $\rho=\gamma/(\gamma+\delta)$. Note that 
$f(U,V)<0$ for all $U,V$ and $\lim_{V\to\infty} f(U,V)=0$ for all $U$. Therefore
(\ref{third}) shows that
\be{aspects}
\begin{array}{lll}
&{\rm (1)} &c_3>0 \mbox{ for } \rho<\tfrac12,\,c_3=0 \mbox{ for } \rho=\tfrac12,\,
c_3<0 \mbox{ for } \rho>\tfrac12,\\
&{\rm (2)} &c_3 \to 0 \mbox{ as } \gamma+\delta \to\infty \mbox{ for fixed } 
\rho\neq\tfrac12 \mbox{ and fixed } \alpha+\beta.
\end{array}
\ee
If $\rho=\tfrac12$, then the dynamics is invariant under swapping the states $0$ and $1$, 
so that $v=0$. If $\rho>\tfrac12$, then $v>0$ for $\alpha-\beta>0$ small enough, but $v$ 
is smaller in the random environment than in the average environment, for which $v=(2\rho-1)
(\alpha-\beta)$ (``slow-down phenomenon''). In the limit $\gamma+\delta\to\infty$ the walk 
sees the average environment.


\subsection{Discussion and outline}
\label{S1.4}

Three classes of models for random walks in dynamic random environments have so 
far been studied in the literature: (1) \emph{space-time random environments}: 
globally updated at each unit of time; (2) \emph{Markovian random environments}: 
independent in space and locally updated according to a single-site Markov chain;
(3) \emph{weak random environments}: small perturbation of homogeneous random walk. 
(See the homepage of Firas Rassoul-Agha [{\sc www.math.utah.edu/$\sim$firas/Research}] 
for an up-to-date list of references.) Our LLN in Theorem~\ref{LLN} is a successful 
attempt to move away from the restrictions. Our expansion of the global speed in 
Theorem~\ref{PERT} is still part of class (3), but it offers some explicit control 
on the coefficients and the domain of convergence of the expansion. 

All papers in the literature deriving LLN's assume an exponential mixing condition
for the dynamic random environment. Cone mixing is one of the weakest mixing 
conditions under which we may expect to be able to derive a LLN via regeneration 
times: no rate of mixing is imposed in (\ref{ConeMix1}). Still, (\ref{ConeMix1}) 
is not optimal because it is a \emph{uniform} mixing condition. For instance, the 
simple symmetric exclusion process, which has a one-parameter family of equilibria 
parameterized by the particle density, is not cone-mixing.       

Both Theorem~\ref{LLN} and \ref{PERT} are easily extended to higher dimension (with
the obvious generalization of cone-mixing), and to random walks whose step rates are
local functions of the environment, i.e., in (\ref{rwtrans}) replace $\xi_t(x)$ by 
$R(\tau_x\xi_t)$, with $\tau_x$ the shift over $x$ and $R$ any cylinder function on 
$\Omega$. It is even possible to allow for steps with a finite range. All that is 
needed is that the total jump rate is independent of the random environment. The 
reader is invited to take a look at the proofs in Sections~\ref{S2} and \ref{S3} to 
see why. In the context of Theorem~\ref{PERT}, the LLN can be extended to a central 
limit theorem (CLT) and to a large deviation principle (LDP), issues which we plan to 
address in future work.


\section{Proof of Theorem~\ref{LLN}}
\label{S2} 

In this section we prove Theorem~\ref{LLN} by adapting the proof of the LLN 
for random walks in \emph{static} random environments developed by Comets and 
Zeitouni~\cite{CoZe04}. The proof proceeds in seven steps. In Section~\ref{S2.1} 
we look at a discrete-time random walk $X$ on $\Z$ in a \emph{dynamic} random 
environment and show that it is equivalent to a discrete-time random walk $Y$ on 
\be{Hdef}
\H=\Z\times\N_0
\ee 
in a \emph{static} random environment that is \emph{directed} in the vertical 
direction. In Section~\ref{S2.2} we show that $Y$ in turn is equivalent to 
a discrete-time random walk $Z$ on $\H$ that suffers \emph{time lapses}, i.e., 
random times intervals during which it does not observe the random environment 
and does not move in the horizontal direction. Because of the cone-mixing property 
of the random environment, these time lapses have the effect of \emph{wiping out 
the memory}. In Section~\ref{S2.3} we introduce \emph{regeneration times} at which, 
roughly speaking, the future of $Z$ becomes independent of its past. Because $Z$ 
is directed, these regeneration times are stopping times. In Section~\ref{S2.4} 
we derive a bound on the moments of the gaps between the regeneration times. 
In Section~\ref{S2.5} we recall a basic coupling property for sequences of random 
variables that are weakly dependent. In Section~\ref{S2.6}, we collect the various 
ingredients and prove the LLN for $Z$, which will immediately imply the LLN for $X$. 
In Section~\ref{S2.7}, finally, we show how the LLN for $X$ can be extended from 
\emph{discrete} time to \emph{continuous} time.

The main ideas in the proof all come from \cite{CoZe04}. In fact, by exploiting 
the directedness we are able to simplify the argument in \cite{CoZe04} considerably. 


\subsection{Space-time embedding}
\label{S2.1} 

Conditional on $\xi$, we define a \emph{discrete-time} random walk on $\Z$
\be{rwdefdis}
X=(X_n)_{n\in\N_0}
\ee 
with transition probabilities 
\be{rwtransdis} 
P^\xi_0\big(X_{n+1}= x+i \mid  X_n=x\big) = \left\{\begin{array}{ll}
p\,\xi_{n+1}(x) + q\,[1-\xi_{n+1}(x)] &\mbox{ if } i=1,\\
q\,\xi_{n+1}(x) + p\,[1-\xi_{n+1}(x)] &\mbox{ if } i=-1,\\
0 &\mbox{ otherwise},
\end{array}
\right. 
\ee 
where $x\in\Z$, $p \in (\tfrac12,1)$, $q=1-p$, and $P^\xi_0$ denotes the law of $X$ starting 
from $X_0=0$ conditional on $\xi$. This is the discrete-time version of the random
walk defined in (\ref{rwdef}--\ref{rwtrans}), with $p$ and $q$ taking over the role
of $\alpha/(\alpha+\beta)$ and $\beta/(\alpha+\beta)$. As in Section~\ref{S1.1}, we 
write $P^\xi_0$ to denote the \emph{quenched} law of $X$ and $\P_{\mu,0}$ to denote 
the \emph{annealed} law of $X$.

Our interacting particle system $\xi$ is assumed to start from an equilibrium 
measure $\mu$ such that the path measure $P^{\mu}$ is stationary and ergodic 
under space-time shifts and is cone-mixing. Given a realization of $\xi$, we 
observe the values of $\xi$ at integer times $n\in\Z$, and introduce a random walk 
on $\H$ 
\be{rwdefY}
Y = (Y_n)_{n\in\N_0}
\ee
with transition probabilities
\be{rw2trans}
P^\xi_{(0,0)}\big(Y_{n+1}= x+e \mid Y_n=x\big)
= \left\{\begin{array}{ll}
p\,\xi_{x_2+1}(x_1) + q\,[1-\xi_{x_2+1}(x_1)]
&\mbox{ if } e=\ell^+,\\
q\,\xi_{x_2+1}(x_1) + p\,[1-\xi_{x_2+1}(x_1)]
&\mbox{ if } e=\ell^-,\\
0 
&\mbox{ otherwise},
\end{array}
\right. 
\ee 
where $x=(x_1,x_2)\in\H$, $\ell^+=(1,1)$, $\ell^-=(-1,1)$, and $P^\xi_{(0,0)}$ denotes 
the law of $Y$ given $Y_0=(0,0)$ conditional on $\xi$. By construction, $Y$ is the random 
walk on $\H$ that moves inside the cone with tip at $(0,0)$ and angle $\tfrac14\pi$, and 
jumps in the directions either $l^+$ or $l^-$, such that
\be{rwequi}
Y_n=(X_n,n), \qquad n\in\N_0.
\ee 
We refer to $P^\xi_{(0,0)}$ as the quenched law of $Y$ and to 
\be{newann} 
\P_{\mu,(0,0)}(\cdot) = \int_{D_\Omega[0,\infty)} P^\xi_{(0,0)}(\cdot)\,P^\mu(\di\xi)
\ee 
as the annealed law of $Y$. If we manage to prove that there exists a $u=(u_1,u_2)\in\R^2$ 
such that
\be{YLLN}
\lim_{n\to\infty} Y_n/n=u \qquad \P_{\mu,(0,0)}-a.s.,
\ee 
then, by (\ref{rwequi}), $u_2=1$, and the LLN in Theorem~\ref{LLN} holds with $v=u_1$.


\subsection{Adding time lapses}
\label{S2.2}

Put $\Lambda=\{(0,0),\ell^+,\ell^-\}$. Let $\epsilon=(\epsilon_i)_{i\in\N}$ be an i.i.d.\ 
sequence of random variables taking values in $\Lambda$ according to the product law 
$W=w^{\otimes\N}$ with marginal 
\be{Qmeas}
w(\epsilon_1=e) = \left\{\begin{array}{ll}
r &\mbox{ if } e\in\{\ell^+,\ell^-\},\\
p &\mbox{ if } e=0,
\end{array}
\right. 
\ee 
with $r=\tfrac12 q$. For fixed $\xi$ and $\epsilon$, introduce a second random walk 
on $\H$
\be{rwdefZ}
Z=(Z_n)_{n\in\N_0}
\ee
with transition probabilities 
\be{EpsWalk}
\begin{aligned}
&\bar{P}_{(0,0)}^{\xi,\epsilon}\big(Z_{n+1}= x+e \mid Z_n=x\big)\\
&\qquad = 1_{\{\epsilon_{n+1}=e\}} + \frac{1}{p}\,1_{\{\epsilon_{n+1}=(0,0)\}}\,
\left[P^\xi_{(0,0)}\big(Y_{n+1}= x+e \mid Y_n=x\big)-r\right],
\end{aligned}
\ee 
where $x\in\H$ and $e\in\{\ell^+,\ell^-\}$, and $\bar{P}_{(0,0)}^{\xi,\epsilon}$ 
denotes the law of $Z$ given $Z_0=(0,0)$ conditional on $\xi,\epsilon$. In words, 
if $\epsilon_{n+1} \in \{\ell^+,\ell^-\}$, then $Z$ takes step $\epsilon_{n+1}$ at
time $n+1$, while if $\epsilon_{n+1}=(0,0)$, then $Z$ copies the step of $Y$.

The quenched and annealed laws of $Z$ defined by
\be{Ybarlaws}
\bar{P}_{(0,0)}^\xi(\cdot) = \int_{\Lambda^\N} \bar{P}_{(0,0)}^{\xi,\epsilon}(\cdot)\,
W(\di\epsilon),
\qquad
\bar{\P}_{\mu,(0,0)}(\cdot) = \int_{D_\Omega[0,\infty)} \bar{P}_{(0,0)}^\xi(\cdot)\,
P^\mu(\di\xi),
\ee
coincide with those of $Y$, i.e.,
\be{Coupling}
\bar{P}_{(0,0)}^\xi(Z \in\,\cdot\,) = P_{(0,0)}^\xi(Y \in\,\cdot\,),
\qquad
\bar{\P}_{\mu,(0,0)}(Z\in\,\cdot\,) =  \P_{\mu,(0,0)}(Y\in\,\cdot\,).
\ee
In words, $Z$ becomes $Y$ when the average over $\epsilon$ is taken. The importance of 
(\ref{Coupling}) is two-fold. First, to prove the LLN for $Y$ in (\ref{YLLN}) it suffices 
to prove the LLN for $Z$. Second, $Z$ suffers time lapses during which its transitions 
are dictated by $\epsilon$ rather than $\xi$. By the cone-mixing property of $\xi$, these 
time lapses will allow $\xi$ to steadily loose memory, which will be a crucial element in 
the proof of the LLN for $Z$.


\subsection{Regeneration times}
\label{S2.3} 

Fix $L\in 2\N$ and define the $L$-vector
\be{LEps}
\epsilon^{(L)}=(\ell^+,\ell^-,\dots,\ell^+,\ell^-),
\ee
where the pair $\ell^+,\ell^-$ is alternated $\tfrac12 L$ times. Given $n\in\N_0$ and 
$\epsilon \in \Lambda^{\N}$ with $(\epsilon_{n+1},\dots,$ $\epsilon_{n+L})=\epsilon^{(L)}$, 
we see from (\ref{EpsWalk}) that (because $\ell^+ +\ell^- = (0,2) = 2\ell$)
\be{InCone}
\bar{P}_{(0,0)}^{\xi,\epsilon}\big(Z_{n+L}= x+L\ell \mid Z_n=x\big)=1,
\qquad x \in \H, 
\ee 
which means that the stretch of walk $Z_n,\dots,Z_{n+L}$ travels in the vertical direction 
$\ell$ irrespective of $\xi$.

Define \emph{regeneration times}
\be{tau1}
\tau_0^{(L)}=0, \qquad 
\tau_{k+1}^{(L)} = \inf\big\{n>\tau_k^{(L)}+L\colon\,
(\epsilon_{n-L},\dots,\epsilon_{n-1})=\epsilon^{(L)}\big\},
\qquad k\in\N.
\ee
Note that these are stopping times w.r.t.\ the filtration $\cG=(\cG_n)_{n\in\N}$
given by
\be{filtration}
\cG_n = \sigma\{\epsilon_i\colon\,1 \leq i\leq n\},\qquad n\in\N.
\ee 
Also note that, by the product structure of $W=w^{\otimes\N}$ defined in (\ref{Qmeas}), 
we have $\tau_k^{(L)}<\infty$ $\bar{\P}_0$-a.s.\ for all $k\in\N$. 

Recall Definition~\ref{PhiMixing} and put
\be{conecor}
\Phi(t) = \sup_{ {A\in\cF_0,\,B\in\cF_t^\theta} \atop {P^\mu(A)>0} }
\Big|P^{\mu}(B\mid A)-P^{\mu}(B)\Big|.
\ee
Cone-mixing is the property that $\lim_{t\to\infty} \Phi(t)=0$ (for all cone angles
$\theta \in (0,\tfrac12\pi)$, in particular, for $\theta=\tfrac14\pi$ needed here). 
Let
\be{FieldHk}
\cH_k=\sigma\Big((\tau_i^{(L)})_{i=0}^k,\,(Z_i)_{i=0}^{\tau_k^{(L)}},\,
(\epsilon_i)_{i=0}^{\tau_k^{(L)}-1},\,\{\xi_t\colon\, 0\leq t\leq \tau_k^{(L)}-L\}\Big), 
\qquad k\in\N.
\ee 
This sequence of sigma-fields allows us to keep track of the walk, the time lapses
and the environment up to each regeneration time. Our main result in the section 
is the following.

\bl{BoundIncr}
For all $L\in 2\N$ and $k\in\N$,
\be{boundInc}
\left\|\,
\bar{\P}_{\mu,(0,0)}\left(Z^{[k]} \in\,\cdot \mid \mathcal{H}_k\right) 
- \bar{\P}_{\mu,(0,0)}\big(Z \in\,\cdot\,\big)\,\right\|_\mathrm{tv} 
\leq \Phi(L),
\ee
where
\be{Ykdef}
Z^{[k]} = \left(Z_{\tau_k^{(L)}+n}-Z_{\tau_k^{(L)}}\right)_{n\in\N_0}
\ee 
and $\|\cdot\|_\mathrm{tv}$ is the total variation norm. 
\el

\begin{proof}
We give the proof for $k=1$. Let $A\in\sigma(\H^{\N_0})$ be arbitrary, and abbreviate 
$1_A=1_{\{Z\in A\}}$. Let $h$ be any $\cH_1$-measurable non-negative random variable. 
Then, for all $x\in\H$ and $n\in\N$, there exists a random variable $h_{x,n}$, measurable 
w.r.t.\ the sigma-field
\be{sigmea}
\sigma\left((Z_i)_{i=0}^n,\,(\epsilon_i)_{i=0}^{n-1},\,
\{\xi_t\colon\,0 \leq t <n-L\}\right),
\ee 
such that $h=h_{x,n}$ on the event $\{Z_n=x,\tau_1^{(L)}=n\}$. Let $E_{P^\mu\otimes W}$
and $\mathrm{Cov}_{P^\mu\otimes W}$ denote expectation and covariance w.r.t.\ $P^\mu\otimes W$, 
and write $\theta_n$ to denote the shift of time over $n$. Then
\be{h}
\begin{aligned}
\bar{\E}_{\mu,(0,0)}\left(h\,\left[1_A \circ \theta_{\tau_1^{(L)}}\right]\right)
&= \sum_{x\in\H, n\in\N} E_{P^\mu\otimes W}\left(
\bar{E}^{\xi,\epsilon}_0\left(h_{x,n}\,[1_A\circ\theta_n]\,
1_{\left\{Z_n=x,\tau_1^{(L)}=n\right\}}\right)\right)\\
&= \sum_{x\in\H, n\in\N} E_{P^\mu\otimes W} \big(f_{x,n}(\xi,\epsilon)\,
g_{x,n}(\xi,\epsilon)\big)\\[0.2cm]
&= \bar{\E}_{\mu,(0,0)}(h)\,\bar{\P}_{\mu,(0,0)}(A) + \rho_A,
\end{aligned}
\ee
where 
\be{fgdefs}
f_{x,n}(\xi,\epsilon) 
= \bar{E}^{\xi,\epsilon}_{(0,0)}\left(h_{x,n}\,1_{\left\{Z_n=x,\tau_1^{(L)}=n\right\}}\right),
\quad 
g_{x,n}(\xi,\epsilon) 
= \bar{P}^{\theta_n\xi,\theta_n\epsilon}_x(A),
\ee
and
\be{rhodef}
\rho_A = \sum_{x\in\H, n\in\N} \mathrm{Cov}_{P^\mu\otimes W}
\big(f_{x,n}(\xi,\epsilon), g_{x,n}(\xi,\epsilon)\big).
\ee
By (\ref{ConeMix1}), we have
\be{Cov}
\begin{aligned}
|\rho_A| &\leq \sum_{x\in\H, n\in\N} 
\big|\mathrm{Cov}_{P^\mu\otimes W}\big(f_{x,n}(\xi,\epsilon), g_{x,n}(\xi,\epsilon)\big)\big|\\
&\leq \sum_{x\in\H, n\in\N} 
\Phi(L)\,E_{P^\mu\otimes W}\big(f_{x,n}(\xi,\epsilon)\big)\,
\sup_{\xi,\epsilon} g_{x,n}(\xi,\epsilon)\\ 
&\leq \Phi(L) \sum_{x\in\H, n\in\N} E_{P^\mu\otimes W}\big(f_{x,n}(\xi,\epsilon)\big)
= \Phi(L)\,\bar{\E}_{\mu,(0,0)}(h).
\end{aligned}
\ee 
Combining (\ref{h}) and (\ref{Cov}), we get
\be{cov3}
\left|\bar{\E}_{\mu,(0,0)}\left(h\,\left[1_A\circ\theta_{\tau_1^{(L)}}\right]\right)
-\bar{\E}_{\mu,(0,0)}(h)\,\bar{\P}_{\mu,(0,0)}(A)\right| \leq \Phi(L)\,\bar{\E}_{\mu,(0,0)}(h). 
\ee 
Now pick $h=1_B$ with $B\in\cH_1$ arbitrary. Then (\ref{cov3}) yields  
\be{add2}
\left|\bar{\P}_{\mu,(0,0)}\left(Z^{[k]} \in A \mid B\right)
-\bar{\P}_{\mu,(0,0)}\left(Z \in A\right)\right| \leq \Phi(L) \mbox{ for all } 
A\in\sigma(\H^{\N_0}),\,B\in\cH_1.
\ee 
There are only countably many cylinders in $\H^{\N_0}$, and so there is a subset 
of $\cH_1$ with $P^\mu$-measure 1 such that, for all $B$ in this set, the above
inequality holds simultaneously for all $A$. Take the supremum over $A$ to get
the claim for $k=1$.

The extension to $k\in\N$ is straightforward. 
\end{proof}


\subsection{Gaps between regeneration times}
\label{S2.4}

Define (recall (\ref{tau1})) 
\be{TimeSeq}
T_k^{(L)} = r^L\left(\tau_k^{(L)}-\tau_{k-1}^{(L)}\right), \qquad k\in\N.
\ee
Note that $T_k^{(L)}$, $k\in\N$, are i.i.d. In this section we prove two lemmas that 
control the moments of these increments.

\bl{IntCon} 
For every $\alpha>1$ there exists an $M(\alpha)<\infty$ such that
\be{IntC}
\sup_{L\in 2\N} \bar{\E}_{\mu,(0,0)}\left([T_1^{(L)}]^{\alpha}\right) \leq M(\alpha). 
\ee
\el

\begin{proof}
Fix $\alpha>1$. Since $T_1^{(L)}$ is independent of $\xi$, we have 
\be{est1}
\bar{\E}_{\mu,(0,0)}\left([T_1^{(L)}]^{\alpha}\right)
= E_W\left([T_1^{(L)}]^{\alpha}\right)
\leq \sup_{L \in 2\N} E_W\left([T_1^{(L)}]^{\alpha}\right),
\ee 
where $E_W$ is expectation w.r.t.\ $W$. Moreover, for all $a>0$, there exists a constant 
$C=C(\alpha,a)$ such that
\be{est2}
[a T_1^{(L)}]^\alpha \leq C\,e^{a T_1^{(L)}},
\ee 
and hence
\be{est3}
\bar{\E}_{\mu,(0,0)}\left([T_1^{(L)}]^{\alpha}\right)
\leq \frac{C}{a^\alpha}\,\sup_{L \in 2\N}\,E_W\left(e^{a T_1^{(L)}}\right).
\ee 
Thus, to get the claim it suffices to show that, for $a$ small enough,
\be{est4}
\sup_{L \in 2\N} \,E_W\left(\erom^{a T_1^{(L)}}\right)<\infty.
\ee 

To prove (\ref{est4}), let
\be{GeomTimes}
I = \inf\big\{m\in\N\colon\,(\epsilon_{mL},\dots,\epsilon_{(m+1)L-1})
=\epsilon^{(L)}\big\}.
\ee
By (\ref{Qmeas}), $I$ is geometrically distributed with parameter $r^L$. Moreover,
$\tau_1^{(L)}\leq (I+1)L$. Therefore 
\be{add3}
\begin{aligned}
E_W\left(\erom^{a T_1^{(L)}}\right)
&= E_W\left(\erom^{a r^L \tau_1^{(L)}}\right)
\leq \erom^{a r^L L}\,E_W\left(\erom^{a r^L I L}\right)\\
&= \erom^{a r^L L} \sum_{j\in\N} (\erom^{a r^L L})^j \,(1-r^L)^{j-1}\,r^L
= \frac{r^L \erom^{2a r^L L}}{\erom^{a r^L L}(1-r^L)},
\end{aligned}
\ee
with the sum convergent for $0<a<(1/r^L L)\log[1/(1-r^L)]$ and tending to zero 
as $L\to\infty$ (because $r<1$). Hence we can choose $a$ small enough so that 
(\ref{est4}) holds.
\end{proof}
 
\bl{lbetaBound}
$\liminf_{L\to\infty} \bar{\E}_{\mu,(0,0)}(T_1^{(L)}) > 0$.
\el

\begin{proof}
Note that $\bar{\E}_{\mu,(0,0)}(T_1^{(L)})<\infty$ by Lemma~\ref{IntCon}. Let 
$N=(N_n)_{n\in\N_0}$ be the Markov chain with state space $S=\{0,1,\dots,L\}$, 
starting from $N_0=0$, such that $N_n=s$ when
\be{add4}
s = 0 \,\vee\, \max\big\{k\in\N\colon\, (\epsilon_{n-k},\dots,\epsilon_{n-1})
= (\epsilon_1^{(L)},\dots,\epsilon_k^{(L)})\big\}
\ee
(with $\max\emptyset=0$). This Markov chain moves up one unit with probability $r$,
drops to $0$ with probability $p+r$ when it is even, and drops to $0$ or $1$ with 
probability $p$, respectively, $r$ when it is odd. Since $\tau_1^{(L)}=\min\{n\in\N_0
\colon\,N_n=L\}$, it follows that $\tau_1^{(L)}$ is bounded from below by a sum of 
independent random variables, each bounded from below by $1$, whose number is 
geometrically distributed with parameter $r^{L-1}$. Hence
\be{bet}
\bar{\P}_{\mu,(0,0)}\left(\tau_1^{(L)} \geq c \,r^{-L}\right) 
\geq (1-r^{L-1})^{\lfloor cr^{-L}\rfloor}.
\ee 
Since
\be{ad1}
\begin{aligned}
\bar{\E}_{\mu,(0,0)}(T_1^{(L)}) &= r^L\,\bar{\E}_{\mu,(0,0)}(\tau_1^{(L)})\\
&\geq r^L\,\bar{\E}_{\mu,(0,0)}\left(\tau_1^{(L)}\,1_{\{\tau_1^{(L)} \geq cr^{-L}\}}\right)
\geq c\,\bar{\P}_{\mu,(0,0)}\left(\tau_1^{(L)} \geq cr^{-L}\right),
\end{aligned}
\ee
it follows that
\be{ad2}
\liminf_{L\rightarrow\infty} \bar{\E}_{\mu,(0,0)}(\tau_1^{(L)}) \geq c\,\erom^{-c/r}.
\ee 
This proves the claim.
\end{proof}


\subsection{A coupling property for random sequences}
\label{S2.5}

In  this section we recall a technical lemma that will be needed in 
Section~\ref{S2.6}. The proof of this lemma is a standard coupling
argument (see e.g.\ Berbee~\cite{Be87}, Lemma 2.1).

\bl{Split}
Let $(U_i)_{i\in\N}$ be a sequence of random variables whose joint probability 
law $P$ is such that, for some marginal probability law $\mu$,
\be{split}
\Big\|P\big(U_i\in\cdot \mid \sigma\{U_j\colon\,1 \leq j < i\}\big)
 - \mu(\cdot)\Big\|_\mathrm{tv} \leq a \quad a.s. \qquad \forall\,i\in\N.
\ee 
Then there exists a sequence of random variables $(\widetilde{U}_i,\Delta_i,
\widehat{U}_i)_{i\in\N}$ satisfying

\medskip
\begin{tabular}{lll}
&(a) &$(\widetilde{U}_i,\Delta_i)_{i\in\N}$ are i.i.d.,\\
&(b) &$\widetilde{U}_i$ has probability law $\mu$,\\
&(c) &$P(\Delta_i=0)=1-a$, $P(\Delta_i=1)=a$,\\
&(d) &$\Delta_i$ is independent of $(\widetilde{U}_j,\Delta_j)_{1 \leq j < i}$
and $\widehat{U}_i$,
\end{tabular}

\medskip\noindent
such that
\be{Xirep} 
U_i = (1-\Delta_i) \widetilde{U}_i + \Delta_i \widehat{U}_i
\quad \mbox{ in distribution}.
\ee
\el


\subsection{LLN for Y}
\label{S2.6} 

Similarly as in (\ref{TimeSeq}), define 
\be{SpaceSeq}
Z_k^{(L)}= r^L \left(Z_{\tau_k^{(L)}}-Z_{\tau_{k-1}^{(L)}}\right),
\qquad k\in\N.
\ee 
In this section we prove the LLN for these increments and this will imply the LLN in 
(\ref{YLLN}).

\begin{proof}
By Lemma~\ref{BoundIncr}, we have
\be{varDis}
\left\|\bar{\P}_{\mu,(0,0)}\big((T_k^{(L)},Z_k^{(L)})
\in\,\cdot\mid\mathcal{H}_{k-1}\big)-\mu^{(L)}(\cdot)\right\|_\mathrm{tv} 
\leq \Phi(L) \quad a.s.\qquad \forall\, k \in \N, 
\ee 
where
\be{ABprop}
\mu^{(L)}(A \times B) = \bar{\P}_{\mu,(0,0)}\big(T_1^{(L)} \in A,Z_1^{(L)}\in B\big)
\qquad \forall\,A \subset r^L\N,\,B \subset r^L\H.
\ee 
Therefore, by Lemma~\ref{Split}, there exists an i.i.d.\ sequence of random variables 
\be{seq1}
(\widetilde{T}_k^{(L)},\widetilde{Z}_k^{(L)},\Delta_k^{(L)})_{k\in\N}
\ee
on $r^L\N \times r^L\H\times\{0,1\}$, where $(\widetilde{T}_k^{(L)},\widetilde{Z}_k^{(L)})$ 
is distributed according to $\mu^{(L)}$ and $\Delta_k^{(L)}$ is Bernoulli distributed with 
parameter $\Phi(L)$, and also a sequence of random variables 
\be{seq2}
(\widehat{Z}_k^{(L)},\widehat{Z}_k^{(L)})_{k\in\N},
\ee 
where $\Delta_k^{(L)}$ is independent of $(\widehat{Z}_k^{(L)},\widehat{Z}_k^{(L)})$ and of
\be{algebraG}
\widetilde{\mathcal{G}}_k
=\sigma\big\{(\widetilde{T}_l^{(L)},\widetilde{Z}_l^{(L)},\Delta_l^{(L)})
\colon\,1 \leq l < k\big\},
\ee
such that
\be{SplitRep}
(T_k^{(L)},Z_k^{(L)})
= (1-\Delta_k^{(L)})\,(\widetilde{T}_k^{(L)},\widetilde{Z}_k^{(L)})
+ \Delta_k^{(L)}\,(\widehat{Z}_k^{(L)},\widehat{Z}_k^{(L)}).
\ee
 
Let
\be{betaMoment}
z_L =\bar{\E}_{\mu,(0,0)}(Z_1^{(L)}),
\ee
which is finite by Lemma~\ref{IntCon} because $|Z_1^{(L)}| \leq T_1^{(L)}$.

\bl{betagamma}
There exists a sequence of numbers $(\delta_L)_{L\in\N_0}$, satisfying $\lim_{L\to\infty}
\delta_L=0$, such that
\be{gammaBound}
\limsup_{n\to\infty}\left|\frac{1}{n} \sum_{k=1}^nZ_k^{(L)}
-z_L \right|<\delta_L \quad \bar{\P}_{\mu,(0,0)}-a.s.
\ee
\el

\begin{proof}
With the help of (\ref{SplitRep}) we can write
\be{TSum}
\frac{1}{n} \sum_{k=1}^n Z_k^{(L)}
=\frac{1}{n} \sum_{k=1}^n \widetilde{Z}_k^{(L)}
-\frac{1}{n} \sum_{k=1}^n \Delta_k^{(L)}\widetilde{Z}_k^{(L)}
+\frac{1}{n} \sum_{k=1}^n \Delta_k^{(L)}\widehat{Z}_k^{(L)}.
\ee
By independence, the first term in the r.h.s.\ of (\ref{TSum}) converges 
$\bar{\P}_{\mu,(0,0)}$-a.s.\ to $z_L$ as $L\to\infty$. H\"older's inequality 
applied to the second term gives, for $\alpha,\alpha'>1$ with 
$\alpha^{-1}+\alpha'^{-1}=1$,
\be{Holder}
\left|\frac{1}{n} \sum_{k=1}^n \Delta_k^{(L)}\widetilde{Z}_k^{(L)}\right|
\leq
\left(\frac{1}{n} \sum_{k=1}^n
\left|\Delta_k^{(L)}\right|^{\alpha'}\right)^{\frac{1}{\alpha'}}
\left(\frac{1}{n} \sum_{k=1}^n
\left|\widetilde{Z}_k^{(L)}\right|^{\alpha}\right)^{\frac{1}{\alpha}}.
\ee
Hence, by Lemma~\ref{IntCon} and the inequality $|\widetilde{Z}_k^{(L)}| 
\leq \widetilde{T}_k^{(L)}$ (compare (\ref{TimeSeq}) and (\ref{SpaceSeq})), 
we have
\be{TSum2}
\limsup_{n\to\infty}\left|\frac{1}{n}
\sum_{k=1}^n \Delta_k^{(L)}\widetilde{Z}_k^{(L)}\right|
\leq \Phi(L)^{\frac{1}{\alpha'}}\,M(\alpha)^{\frac{1}{\alpha}}
\quad \bar{\P}_{\mu,(0,0)}-a.s.
\ee

It remains to analyze the third term in the r.h.s.\ of (\ref{TSum}). Since $|\Delta_k^{(L)}
\widehat{Z}_k^{(L)}|\leq Z_k^{(L)}$, it follows from Lemma~\ref{IntCon} that 
\be{MBound}
\begin{aligned}
M(\alpha) &\geq \bar{\E}_{\mu,(0,0)}\big(|Z_k^{(L)}|^\alpha\big)\\
&\geq \bar{\E}_{\mu,(0,0)}\left(|\Delta_k^{(L)}\widehat{Z}_k^{(L)}|^{\alpha}\mid
\widetilde{\mathcal{G}}_k\right)
= \Phi(L)\,\,\bar{\E}_{\mu,(0,0)}\left(|\widehat{Z}_k^{(L)}|^{\alpha}\mid
\widetilde{\mathcal{G}}_k\right) \quad a.s.
\end{aligned}
\ee 
Next, put $\widehat{Z}_k^{*(L)}=\bar{\E}_{\mu,(0,0)}(\widehat{Z}_k^{(L)}\mid 
\widetilde{\mathcal{G}}_k)$ and note that 
\be{martingale}
M_n = \frac{1}{n} \sum_{k=1}^n \Delta_k^{(L)} 
\left(\widehat{Z}_k^{(L)}-\widehat{Z}_k^{*(L)}\right)
\ee
is a mean-zero martingale w.r.t.\ the filtration $\widetilde{\mathcal{G}}=(\widetilde
{\mathcal{G}}_k)_{k\in\N}$. By the Burkholder-Gundy maximal inequality (Williams~\cite{Wi91},
(14.18)), it follows that, for $\beta =\alpha \wedge 2$,
\be{martineq}
\begin{aligned}
&\bar{\E}_{\mu,(0,0)}\left(\Big|\sup_{n\in\N} M_n\Big|^\beta\right)
\leq C(\beta)\,\bar{\E}_{\mu,(0,0)}\Big(\sum_{k\in\N}
\frac{[\Delta_k^{(L)}(\widehat{Z}_k^{(L)}-\widehat{Z}_k^{*(L)})]^2}{k^2}\Big)^{\beta/2}\\
&\leq C(\beta) \sum_{k\in\N} \bar{\E}_{\mu,(0,0)}
\left(\frac{|\Delta_k^{(L)}(\widehat{Z}_k^{(L)}-\widehat{Z}_k^{*(L)})|^\beta}{k^\beta}\right)
\leq C'(\beta), 
\end{aligned}
\ee 
for some constants $C(\beta),C'(\beta)<\infty$. Hence $M_n$ a.s.\ converges to an integrable 
random variable as $n\to\infty$, and by Kronecker's lemma $\lim_{n\to\infty} M_n=0$ a.s. 
Moreover, if $\Phi(L)>0$, then by Jensen's inequality and (\ref{MBound}) we have
\be{add5}
|\widehat{Z}_k^{*(L)}| \leq
\left[\bar{\E}_{\mu,(0,0)}\left(\big|\widehat{Z}_k^{(L)}\big|^{\alpha} \mid
\widetilde{\mathcal{G}}_k\right)\right]^{\frac{1}{\alpha}}
\leq \left(\frac{M(\alpha)}{\Phi(L)}\right)^{\frac{1}{\alpha}}
\quad \bar{\P}_{\mu,(0,0)}-a.s.
\ee
Hence
\be{TSum3}
\left|\frac{1}{n} \sum_{k=1}^n \Delta_k^{(L)}\widehat{Z}_k^{*(L)}\right|
\leq \left(\frac{M(\alpha)}{\Phi(L)}\right)^{\frac{1}{\alpha}}
\frac{1}{n} \sum_{k=1}^n \Delta_k^{(L)}.
\ee 
As $n\to\infty$, the r.h.s.\ converges $\bar{\P}_{\mu,(0,0)}$-a.s.\ to $M(\alpha)^{\frac{1}
{\alpha}} \Phi(L)^{\frac{1}{\alpha'}}$. Therefore, recalling (\ref{TSum3}) and choosing 
$\delta_L=2M(\alpha)^{\frac{1}{\alpha}}\Phi(L)^{\frac{1}{\alpha'}}$, we get the claim.
\end{proof}

Finally, since $\widetilde{Z}_k^{(L)} \geq r^L$ and 
\be{LLNTs}
\frac1n \sum_{k=1}^n T_k^{(L)} = t_L = \bar{\E}_{\mu,(0,0)}(T_1^{(L)}) >0
\quad \bar{\P}_{\mu,(0,0)}-a.s., 
\ee
Lemma~\ref{betagamma} yields
\be{C1}
\limsup_{n\to\infty} \left|\frac{\frac{1}{n}
\sum_{k=1}^n Z_k^{(L)}} {\frac{1}{n} \sum_{k=1}^n
Z_k^{(L)}}-\frac{z_L}{t_L}\right|<C_1\,\delta_L 
\quad \bar{\P}_{\mu,(0,0)}-a.s.
\ee 
for some constant $C_1<\infty$ and $L$ large enough. By (\ref{TimeSeq}) and (\ref{SpaceSeq}),
the quotient of sums in the l.h.s.\ equals $Z_{\tau_n^{(L)}}/\tau_n^{(L)}$. It therefore 
follows from a standard interpolation argument that
\be{C2}
\limsup_{n\to\infty}\left|\frac{Z_n}{n}-\frac{z_L}{t_L}\right|
<C_2\,\delta_L \quad \bar{\P}_{\mu,(0,0)}-a.s.
\ee 
for some constant $C_2<\infty$ and $L$ large enough. This implies the existence of 
the limit $\lim_{L\to\infty} z_L/t_L$, as well as the fact that $\lim_{n\to\infty} 
Z_n/n = u$ $\bar{\P}_{\mu,(0,0)}$-a.s., which in view of (\ref{Coupling}) is equivalent 
to the statement in (\ref{YLLN}) with $u=(v,1)$.
\end{proof}


\subsection{From discrete to continous time}
\label{S2.7}

It remains to show that the LLN derived in Sections~\ref{S2.1}--\ref{S2.6} for
the discrete-time random walk defined in (\ref{rwdefdis}--\ref{rwtransdis}) can 
be extended to the continuous-time random walk defined in (\ref{rwdef}--\ref{rwtrans}).

Let $\chi=(\chi_n)_{n\in\N_0}$ denote the jump times of the continuous-time random 
walk $X=(X_t)_{t \geq 0}$ (with $\chi_0=0$). Let $Q$ denote the law of $\chi$. The 
increments of $\chi$ are i.i.d.\ random variables, independent of $\xi$, whose 
distribution is exponential with mean $1/(\alpha+\beta)$. Define
\be{tildexiX}
\begin{array}{lllll}
\xi^* &=& (\xi^*_n)_{n\in\N_0} \quad \mbox{ with } \quad \xi^*_n &=& \xi_{\chi_n},\\ 
X^* &=& (X^*_n)_{n\in\N_0} \quad \mbox{ with} \quad X^*_n &=& X_{\chi_n}. 
\end{array}
\ee
Then $X^*$ is a discrete-time random walk in a discrete-time random environment of
the type considered in Sections~\ref{S2.1}--\ref{S2.6}, with $p=\alpha/(\alpha+\beta)$
and $q=\beta/(\alpha+\beta)$. Lemma~\ref{conediscont} below shows that the cone-mixing 
property of $\xi$ carries over to $\xi^*$ under the joint law $P^\mu \times Q$. Therefore 
we have (recall (\ref{annealed}))
\be{LLNX*}
\lim_{n\to\infty} X^*_n/n = v^* \quad \mbox{ exists } (\P_{\mu,0} \times Q)-a.s.
\ee 
Since $\lim_{n\to\infty} \chi_n/n = 1/(\alpha+\beta)$ $Q$-a.s., it follows that
\be{LLNX*ext}
\lim_{n\to\infty} X_{\chi_n}/\chi_n = (\alpha+\beta)v^* \quad \mbox{ exists } 
(\P_{\mu,0} \times Q)-a.s.
\ee
A standard interpolation argument now yields (\ref{LLN1}) with $v=(\alpha+\beta)v^*$.

\bl{conediscont}
If $\xi$ is cone-mixing with angle $\theta>\arctan(\alpha+\beta)$, then $\xi^*$ is 
cone-mixing with angle $\tfrac14\pi$.
\el

\begin{proof}
Fix $\theta>\arctan(\alpha+\beta)$, and put $c=c(\theta) = \cot\theta < 1/(\alpha+\beta)$.
Recall from (\ref{Cone}) that $C_t^\theta$ is the cone with angle $\theta$ whose tip is
at $(0,t)$. For $M\in\N$, let $C_{t,M}^\theta$ be the cone obtained from $C_t^\theta$ by 
extending the tip to a rectangle with base $M$, i.e.,
\be{CtNdef}
C_{t,M}^\theta = C_t^\theta \cup \{([-M,M] \cap \Z) \times [t,\infty)\}.
\ee
Because $\xi$ is cone-mixing with angle $\theta$, and
\be{coneinlc}
C_{t,M}^\theta \subset C_{t-cM}^\theta, \qquad M\in\N,
\ee 
$\xi$ is cone-mixing with angle $\theta$ and base $M$, i.e., (\ref{ConeMix1}) holds with 
$C_t^\theta$ replaced by $C_{t,M}^\theta$. This is true for every $M\in\N$.

Define, for $t \geq 0$ and $M\in\N$,
\be{sigalgdef1}
\begin{aligned}
\cF_t^\theta &= \sigma\big\{\xi_s(x)\colon\,(x,s) \in C_t^\theta\big\},\\
\cF_{t,M}^\theta &= \sigma\big\{\xi_s(x)\colon\,(x,s) \in C_{t,M}^\theta\big\},\\
\end{aligned}
\ee
and, for $n\in\N$,
\be{sigalgdef2}
\begin{aligned}
\cF^*_n &= \sigma\big\{\xi_m^*(x)\colon\,(x,m) \in C_n^{\tfrac14\pi}\big\},\\
\cG_n &= \sigma\big\{\chi_m\colon\,m \geq n\big\},  
\end{aligned}
\ee
where $C_n^{\tfrac14\pi}$ is the discrete-time cone with tip $(0,n)$ and angle $\tfrac14\pi$. 

Fix $\delta>0$. Then there exists an $M=M(\delta)\in\N$ such that $Q(D[M]) \geq 1-\delta$
with $D[M] = \{\chi_n/n \geq c \,\,\forall\,n \geq M\}$. For $n\in\N$, define
\be{Dndef}
D_n = \big\{\chi_n/n \geq c\big\} \cap \sigma^n D[M],
\ee
where $\sigma$ is the left-shift acting on $\chi$. Since $c<1/(\alpha+\beta)$, we have
$P(\chi_n/n \geq c) \geq 1-\delta$ for $n \geq N = N(\delta)$, and hence $P(D_n) \geq 
(1-\delta)^2 \geq 1-2\delta$ for $n \geq N = N(\delta)$,. Next, obeserve that
\be{key}
B \in \cF^*_n \Longrightarrow B \cap D_n \in \cF_{cn,M}^\theta \otimes \cG_n
\ee
(the r.h.s.\ is the product sigma-algebra). Indeed, on the event $D_n$ we have $\chi_m \geq cm$
for $m \geq n+M$, which implies that, for $m \geq M$,
\be{coneemb}
(x,m) \in C_n^{\tfrac14\pi} \Longrightarrow |x|+m \geq n \Longrightarrow c|x|+\chi_n \geq cn
\Longrightarrow (x,\chi_m) \in C_{cn,M}^\theta.
\ee

Now put $\bar{P}^\mu = P^\mu \otimes Q$ and, for $A \in \cF_0$ with $P^\mu(A)>0$ and 
$B \in \cF^*_n$ estimate
\be{Pthree}
|\bar{P}^\mu(B \mid A) - \bar{P}^\mu(B)| \leq I + II + III
\ee
with
\be{three}
\begin{aligned}
I &= |\bar{P}^\mu(B \mid A) - \bar{P}^\mu(B \cap D_n \mid A)|,\\
II &= |\bar{P}^\mu(B \cap D_n \mid A) - \bar{P}^\mu(B \cap D_n)|,\\
III &= |\bar{P}^\mu(B \cap D_n) - \bar{P}^\mu(B)|.
\end{aligned}
\ee
Since $D_n$ is independent of $A,B$ and $P(D_n) \geq 1-2\delta$, it follows that $I \leq 2\delta$ 
and $III \leq 2\delta$ uniformly in $A$ and $B$. To bound $II$, we use (\ref{key}) to estimate
\be{IIest}
II \leq \sup_{ {A \in \cF_0,\,B' \in \cF_{cn,M}^\theta \otimes \cG_n} \atop {P^\mu(A)>0} }
|\bar{P}^\mu(B' \mid A) - \bar{P}^\mu(B')|.
\ee
But the r.h.s.\ is bounded from above by
\be{contrbd}
\sup_{ {A \in \cF_0,\,B'' \in \cF_{cn,M}^\theta} \atop {P^\mu(A)>0} }
|P^\mu(B'' \mid A) - P^\mu(B'')|
\ee
because, for every $B'' \in \cF_{cn,M}^\theta$ and $C \in \cG_n$, 
\be{estPs}
|\bar{P}^\mu(B'' \times C \mid A) - \bar{P}^\mu(B'' \times C)| 
= |[P^\mu(B'' \mid A) - P^\mu(B'')]\,Q(C)| \leq |P^\mu(B''\mid A)-P^\mu(B'')|,
\ee
where we use that $C$ is independent of $A,B''$.

Finally, because $\xi$ is cone-mixing with angle $\theta$ and base $M$, (\ref{contrbd}) tends
to zero as $n\to\infty$, and so by combining (\ref{Pthree}--\ref{contrbd}) we get
\be{}
\limsup_{n\to\infty} \sup_{ {A \in \cF_0,\,B \in \cF^*_n} \atop {P^\mu(A)>0} } 
|\bar{P}^\mu(B \mid A) - \bar{P}^\mu(B)| \leq 4\delta.
\ee
Now let $\delta \downarrow 0$ to obtain that $\xi^*$ is cone mixing with angle $\tfrac14\pi$.
\end{proof}


\section{Series expansion for $M<\epsilon$}
\label{S3}

Throughout this section we assume that the dynamic random environment $\xi$ falls 
in the regime for which $M<\epsilon$ (recall (\ref{cTdef}). In Section~\ref{S3.1} we 
define the \emph{environment process}, i.e., the environment as seen relative to the 
position of the random walk. In Section~\ref{S3.2} we prove that this environment 
process has a \emph{unique ergodic equilibrium} $\mu_e$, and we derive a series 
expansion for $\mu_e$ in powers of $\alpha-\beta$ that converges when $\alpha-\beta
<\tfrac12(\epsilon-M)$. In Section~\ref{S3.3} we use the latter to derive a series 
expansion for the global speed $v$ of the random walk.


\subsection{Definition of the environment process}
\label{S3.1}

Let $X=(X_t)_{t \geq 0}$ be the random walk defined in (\ref{rwdef}--\ref{rwtrans}).
For $x\in\Z$, let $\tau_x$ denote the shift of space over $x$.

\bd{EnvPRo}
The environment process is the Markov process $\zeta=(\zeta_t)_{t \geq 0}$ with
state space $\Omega$ given by
\be{envpro}
\zeta_t = \tau_{X_t}\xi_t, \qquad t \geq 0, 
\ee
where
\be{tau}
(\tau_{X_t}\xi_t)(x) = \xi_t(x+X_t), \qquad x \in \Z,\,t \geq 0.
\ee
Equivalently, if $\xi$ has generator $L_\mathrm{IPS}$, then $\zeta$ has
generator $L$ given by
\be{GenEnvPro}
(Lf)(\eta) =
c^+(\eta)\big[f(\tau_1\eta)-f(\eta)\big]
+c^-(\eta)\big[f(\tau_{-1}\eta)-f(\eta)\big] +
(L_\mathrm{IPS}f)(\eta), \quad \eta \in \Omega,
\ee
where $f$ is an arbitrary cylinder function on $\Omega$ and
\be{dxrate}
\begin{aligned}
c^+(\eta) &= \alpha\,\eta(0)+\beta\,[1-\eta(0)],\\
c^-(\eta) &= \beta\,\eta(0)+\alpha\,[1-\eta(0)].
\end{aligned}
\ee
\ed

Let $S=(S(t))_{t\geq 0}$ be the semigroup associated with the generator $L$. Suppose
that we manage to prove that $\zeta$ is ergodic, i.e., there exists a unique probability
measure $\mu_e$ on $\Omega$ such that, for any cylinder function $f$ on $\Omega$,
\be{ergodicMeas}
\lim_{t\to\infty} (S(t)f)(\eta) = \langle
f\rangle_{\mu_e} \qquad \forall\,\eta\in\Omega,
\ee
where $\langle \cdot \rangle_{\mu_e}$ denotes expectation w.r.t. $\mu_e$. Then, picking
$f = \phi_0$ with $\phi_0(\eta)=\eta(0)$, $\eta\in\Omega$, we have
\be{envDensity} 
\lim_{t\to\infty} (S(t)\phi_0)(\eta) = \langle
\phi_0\rangle_{\mu_e} = \widetilde{\rho} \qquad
\forall\,\eta\in\Omega
\ee
for some $\widetilde{\rho} \in [0,1]$, which represents the limiting probability that 
$X$ is on an occupied site given that $\xi_0=\zeta_0=\eta$ (note that $(S(t)\phi_0)(\eta) 
= E^\eta(\zeta_t(0))=E^\eta(\xi_t(X_t))$). 

Next, let $N_t^+$ and $N_t^-$ be the number of shifts to the right, respectively, left up 
to time $t$ in the environment process. Then $X_t=N_t^+-N_t^-$. Since $M_t^j = N_t^j
-\int_0^t c^j(\eta_s)\,\di s$, $j\in\{+,-\}$, are martingales with stationary and ergodic 
increments, we have
\be{Martingale}
X_t = M_t+(\alpha-\beta) \int_0^t \big(2\eta_s(0)-1\big)\,\di s
\ee 
with $M_t=M_t^+-M_t^-$ a martingle with stationary and ergodic increments. It 
follows from (\ref{envDensity}--\ref{Martingale}) that
\be{LLNlim} 
\lim_{t\to\infty} X_t/t = (2\widetilde{\rho}-1)(\alpha-\beta) \qquad \mu-a.s. 
\ee

In Section~\ref{S3.2} we prove the existence of $\mu_e$, and show that it can be expanded 
in powers of $\alpha-\beta$ when $\alpha-\beta<\tfrac12(\epsilon-M)$. In Section~\ref{S3.3} 
we use this expansion to obtain an expansion of $\widetilde\rho$.


\subsection{Unique ergodic equilibrium measure for the environment process}
\label{S3.2}

In Section~\ref{S3.2.1} we prove four lemmas controlling the evolution of $\zeta$. In 
Section~\ref{S3.2.2} we use these lemmas to show that $\zeta$ has a unique ergodic 
equilibrium measure $\mu_e$ that can be expanded in powers of $\alpha-\beta$, provided 
$\alpha-\beta<\tfrac12(\epsilon-M)$.

We need some notation. Let $\|\cdot\|_\infty$ be the sup-norm on $C(\Omega)$. Let 
$\trn\cdot\trn$ be the triple norm on $\Omega$ defined as follows. For $x\in\Z$ and 
a cylinder function $f$ on $\Omega$, let 
\be{norm}
\Delta_f(x) = \sup_{\eta\in\Omega} |f(\eta^x)-f(\eta)|
\ee
be the maximum variation of $f$ at $x$, where $\eta^x$ is the configuration obtained 
from $\eta$ by flipping the state at site $x$, and put
\be{TripleNorm}
\trn f \trn =\sum_{x\in\Z}\Delta_f(x).
\ee
It is easy to check that, for arbitrary cylinder functions $f$ and $g$ on $\Omega$,
\be{TripleProduct}
\trn f g\trn\leq \|f\|_\infty\, \trn g \trn + \|g\|_\infty\, \trn f \trn.
\ee


\subsubsection{Decomposition of the generator of the environment process}
\label{S3.2.1}

\bl{GenDecomposition} 
Assume {\rm (\ref{Pmuass})} and suppose that $M<\epsilon$. Write the generator of the 
environment process $\zeta$ defined in {\rm (\ref{GenEnvPro})} as 
\be{RevIrr} 
L=L_0+L_* = (L_\mathrm{SRW} + L_\mathrm{IPS}) + L_* , 
\ee 
where 
\be{RevIrrev}
\begin{aligned}
(L_\mathrm{SRW} f)(\eta) &= \tfrac12(\alpha+\beta)\,
\Big[f(\tau_1\eta)+f(\tau_{-1}\eta)-2f(\eta)\Big],\\
(L_* f)(\eta) &= \tfrac12(\alpha-\beta)\,
\Big[f(\tau_1\eta)-f(\tau_{-1}\eta)\Big]\,\big(2\eta(0)-1\big).
\end{aligned}
\ee 
Then $L_0$ is the generator of a Markov process that still has $\mu$ as an equilibrium, 
and that satisfies 
\be{TrnBound} 
\trn S_0(t)f \trn  \leq \erom^{-ct}\,\trn f \trn 
\ee 
and
\be{SupnBound} 
\|S_0(t)f-\langle f \rangle_\mu \|_\infty \leq C\,\erom^{-ct}\,\trn f \trn, 
\ee 
where $S_0=(S_0(t))_{t \geq 0}$ is the semigroup associated with the generator $L_0$,
$c=\epsilon-M$, and $C<\infty$ is a positive constant. 
\el

\begin{proof}
Note that $L_\mathrm{SRW}$ and $L_\mathrm{IPS}$ commute. Therefore, for an arbitrary 
cylinder function $f$ on $\Omega$, we have 
\be{ExpConv1} 
\trn S_0(t)f \trn = \trn \erom^{t L_\mathrm{SRW}}
\big(\erom^{t L_\mathrm{IPS}}f\big) \trn 
\leq \trn \erom^{t L_\mathrm{IPS}} f \trn  \leq \erom^{-ct}\,\trn f \trn, 
\ee 
where the first inequality uses that $\erom^{t L_\mathrm{SRW}}$ is a contraction 
semigroup, and the second inequality follows from the fact that $\xi$ falls in 
the regime $M<\epsilon$ (see Liggett~\cite{Li85}, Theorem I.3.9). The inequality 
in (\ref{SupnBound}) follows by a similar argument. Indeed, 
\be{ExpConv2} \|S_0(t)f-\langle f \rangle_\mu \|_\infty 
= \|\erom^{t L_\mathrm{SRW}}\big(\erom^{t L_\mathrm{IPS}}f\big)
-\langle f \rangle_\mu \|_\infty \leq\|\erom^{tL_\mathrm{IPS}}f
-\langle f \rangle_\mu \|_\infty 
\leq C\,\erom^{-ct}\,\trn f \trn, 
\ee 
where the last inequality again uses that $\xi$ falls in the regime $M<\epsilon$ 
(see Liggett \cite{Li85}, Theorem I.4.1). The fact that $\mu$ is an equilibrium measure 
is trivial, since $L_\mathrm{SRW}$ only acts on $\eta$ by shifting it.
\end{proof}

Note that $L_\mathrm{SRW}$ is the generator of simple random walk on $\Z$ jumping at rate 
$\alpha+\beta$. We view $L_0$ as the generator of an unperturbed Markov process and $L_*$ 
as a perturbation of $L_0$. The following lemma gives us control of the latter. 

\bl{Bounds}
For any cylinder function $f$ on $\Omega$,
\be{IrrSupBound} 
\|L_* \,f\|_\infty \leq (\alpha-\beta)\|f\|_\infty
\ee
and
\be{IrrTripleBound}
\trn L_*\,f \trn \leq 2(\alpha-\beta)\,\trn f \trn 
\quad \text{ if } \langle f \rangle_\mu=0.
\ee
\el

\begin{proof}
To prove (\ref{IrrSupBound}), estimate 
\be{sn}
\begin{aligned}
\|L_*\,f\|_\infty &= \tfrac12(\alpha-\beta)\,
\| \big[f(\tau_1\,\cdot)-f(\tau_{-1}\,\cdot)\big]\,\big(2\phi_0(\cdot)-1\big)\|_\infty\\
&\leq \tfrac12(\alpha-\beta)\, \|f(\tau_1\,\cdot)+ f(\tau_{-1}\,\cdot)\|_\infty
\leq (\alpha-\beta)\,\|f\|_\infty.
\end{aligned}
\ee
To prove (\ref{IrrTripleBound}), recall (\ref{RevIrrev}) and estimate
\be{ppp1}
\begin{aligned}
\trn L_* f \trn &= \tfrac12(\alpha-\beta)\,
\trn \big[f(\tau_1\,\cdot)-f(\tau_{-1}\,\cdot)\big]\,\big(2\phi_0(\cdot)-1\big)\trn\\
&\leq \tfrac12(\alpha-\beta)\, \Big\{\trn f(\tau_1\cdot)(2\phi_0(\cdot)-1)\trn
+\trn f(\tau_{-1}\cdot)(2\phi_0(\cdot)-1)\trn\Big\}\\
&\leq (\alpha-\beta)\,\Big(\|f\|_\infty\,\trn\,(2\phi_0-1)\trn
+\trn f \trn\,\|(2\phi_0-1)\|_\infty\Big)\\
&= (\alpha-\beta)\,\Big(\|f\|_\infty +\trn f \trn\Big)
\leq 2(\alpha-\beta)\trn f \trn,
\end{aligned}
\ee
where the second inequality uses (\ref{TripleProduct}) and the third inequality follows 
from the fact that $\|f\|_\infty \leq \trn f \trn$ for any $f$ such that $\langle f 
\rangle_\mu=0$.
\end{proof}

We are now ready to expand the semigroup $S$ of $\zeta$. Henceforth abbreviate
\be{cdef}
c = \epsilon-M.
\ee

\bl{Expansion}
Let $S_0=(S_0(t))_{t \geq 0}$ be the semigroup associated with the generator $L_0$
defined in {\rm (\ref{RevIrrev})}. Then, for any $t \geq 0$ and any cylinder function
$f$ on $\Omega$,
\be{expansion}
S(t)f = \sum_{n\in\N} g_n(t,f),
\ee
where
\be{g_n}
g_1(t,f) = S_0(t)f \quad
\text{and}\quad g_{n+1}(t,f) = \int_0^t
S_0(t-s)\,L_*\,g_n(s,f)\,\di s, \qquad n\in\N.
\ee
Moreover, for all $n\in\N$,
\be{gSupBound}
\|g_n(t,f)\|_\infty \leq \trn f \trn\,
\Big(\frac{2(\alpha-\beta)}{c}\Big)^{n-1}
\ee
and
\be{gTripleBound}
\trn g_n(t,f) \trn \leq
\,\erom^{-ct}\,\frac{[2(\alpha-\beta)t]^{n-1}}{(n-1)!}\, \trn f \trn,
\ee
where $0!=1$. In particular, for all $t>0$ and $\alpha-\beta<\tfrac12 c$ 
the series in {\rm (\ref{expansion})} converges uniformly in $\eta$.
\el

\begin{proof}
Since $L=L_0+L_*$, Dyson's formula gives
\be{Dyson1} 
\erom^{tL}\,f=\erom^{tL_0}f + \int_0^t
\erom^{(t-s)L_0}\,L_*\,\erom^{sL}\,f\,\di s,
\ee
which, in terms of semigroups, reads
\be{Dyson2}
S(t)f = S_0(t)f + \int_0^t S_0(t-s)L_*\,S(s)f\,\di s.
\ee
The expansion in (\ref{expansion}--\ref{g_n}) follows from (\ref{Dyson2})
by induction on $n$. 

We next prove (\ref{gTripleBound}) by induction on $n$. For $n=1$ the 
claim is immediate. Indeed, by Lemma~\ref{GenDecomposition} we have the 
exponential bound
\be{n=1}
\trn g_1(t,f) \trn = \trn S_0(t)f \trn \leq \erom^{-ct}\,\trn f \trn.
\ee
Suppose that the statement in (\ref{gTripleBound}) is true up to $n$. Then
\be{induction1}
\begin{aligned}
\trn g_{n+1}(t,f) \trn 
&=\trn \int_0^t S_0(t-s)\,L_*\,g_n(s,f)\,\di s\,\trn\\
&\leq \int_0^t \trn S_0(t-s)\,L_*\,g_n(s,f)\trn\,\di s\\
&\leq \int_0^t \erom^{-c(t-s)}\,\trn L_*\,g_n(s,f)\trn\,\di s\\
&= \int_0^t \erom^{-c(t-s)}\,\trn L_*\,\big(g_n(s,f)
-\langle g_n(s,f)\rangle_\mu\big) \trn\, \di s\\
&\leq 2(\alpha-\beta)\,\int_0^t \erom^{-c(t-s)}\,\trn g_n(s,f)\trn\,\di s,\\
&\leq \trn f \trn\,\erom^{-ct}\, [2(\alpha-\beta)]^n\,
\int_0^t \frac{s^{n-1}}{(n-1)!}\,\di s\\
&= \trn f \trn\,\erom^{-ct}\,\frac{[2(\alpha-\beta)t]^n}{n!},
\end{aligned}
\ee
where the third inequality uses (\ref{IrrTripleBound}), and the fourth inequality 
relies on the induction hypothesis. 

Using (\ref{gTripleBound}), we can now prove (\ref{gSupBound}). Estimate 
\be{sup}
\begin{aligned}
\|g_{n+1}(t,f) \|_\infty 
&=\left\| \int_0^t S_0(t-s)\,L_*\,g_n(s,f) \di s\right\|_\infty\\ 
&\leq \int_0^t \,\|L_*\,g_n(s,f)\|_\infty\,\di s\\
&= \int_0^t \,\big\| L_*\,\big(g_n(s,f)-\langle g_n(s)\rangle_\mu\big)\big\|_\infty\,\di s\\ 
&\leq (\alpha-\beta) \int_0^t \big\| g_n(s,f)-\langle g_n(s,f)\rangle_\mu\big\|_\infty\,\di s\\
&\leq (\alpha-\beta) \int_0^t \,\trn g_n(s,f)\trn\,\di s\\ 
&\leq (\alpha-\beta)\trn f \trn \int_0^t \erom^{-cs}\,
\frac{[2(\alpha-\beta)s]^{n-1}}{(n-1)!}\,\di s\\
&\leq \trn f \trn\Big(\frac{2(\alpha-\beta)}{c}\Big)^n,
\end{aligned}
\ee
where the first inequality uses that $S_0(t)$ is a contraction semigroup, while the 
second and fourth inequality rely on (\ref{IrrSupBound}) and (\ref{gTripleBound}).
\end{proof}

We next show that the functions in (\ref{expansion}) are uniformly close to their average 
value.

\bl{h_nBound}
Let
\be{h_n}
h_n(t,f) = g_n(t,f)-\langle g_n(t,f) \rangle_\mu, \quad t \geq
0,\,n\in\N.
\ee
Then 
\be{hBound}
\| h_n(t,f) \|_\infty \leq
C\,\erom^{-ct}\,\frac{[2(\alpha-\beta)t]^{n-1}}{(n-1)!}\, \trn f
\trn,
\ee
for some $C<\infty$ ($0!=1$).
\el

\begin{proof}
Note that $\trn h_n(t,f)\trn = \trn g_n(t,f) \trn$ for  $t\geq 0$ and
$n\in\N$, and estimate
\be{induction2}
\begin{aligned}
\|h_{n+1}(t,f) \|_\infty &=\left\| \int_0^t
\Big(S_0(t-s)\,L_*\,g_n(s,f)
-\langle L_*\,g_n(s,f) \rangle_\mu \Big)\,\di s\right\|_\infty\\
&\leq C \int_0^t \erom^{-c(t-s)}\,\trn L_*\,g_n(s,f)\trn\,\di s\\
&= C \int_0^t \erom^{-c(t-s)}\,\trn L_*\,h_n(s,f)\trn\,\di s\\
&\leq C\,2(\alpha-\beta) \int_0^t \erom^{-c(t-s)}\,\trn h_n(s,f)\trn\,\di s\\
&\leq C\,\trn f \trn\,\erom^{-ct}\,[2(\alpha-\beta)]^n
\int_0^t \frac{s^{n-1}}{(n-1)!}\,\di s\\
&= C\,\trn f \trn\,\erom^{-ct}\,\frac{[2(\alpha-\beta)t]^n}{n!},
\end{aligned}
\ee
where the first inequality uses (\ref{SupnBound}), while the second and third inequality 
rely on (\ref{IrrTripleBound}) and (\ref{gTripleBound}).
\end{proof}


\subsubsection{Expansion of the equilibrium measure of the environment process}
\label{S3.2.2}

We are finally ready to state the main result of this section.

\bt{ErgodicMeasure}
For $\alpha-\beta<\tfrac12 c$, the environment process $\zeta$ has a unique invariant 
measure $\mu_e$. In particular, for any cylinder function $f$ on $\Omega$,
\be{ErgMeas}
\langle f \rangle_{\mu_e} = \lim_{t\to\infty} \langle
S(t)f \rangle_\mu = \sum_{n\in\N} \lim_{t\to\infty}\langle g_n(t,f) \rangle_\mu.
\ee
\et

\begin{proof}
By Lemma~\ref{h_nBound}, we have
\be{convergence}
\begin{aligned}
\left\| S(t)f - \left\langle S(t)f \right\rangle_\mu\right\|_\infty
&=\left\|\sum_{n\in\N} g_{n}(t,f)-  \langle \sum_{n\in\N}g_n(t,f)
\rangle_\mu\right\|_\infty= \left\|\sum_{n\in\N}
h_{n}(t,f)\right\|_\infty
\\&\leq \sum_{n\in\N} \|h_{n}(t,f)\|_\infty
\leq C\,\erom^{-ct}\,\trn f \trn
\sum_{n\in\N}\frac{[2(\alpha-\beta)t]^n}{n!} \\
&= C \trn f \trn\,
\erom^{-t[c-2(\alpha-\beta)]}.
\end{aligned}
\ee
Since $\alpha-\beta<\tfrac12 c$, we see that the r.h.s.\  of (\ref{convergence}) tends to zero 
as $t\to\infty$. Consequently, the l.h.s.\ tends to zero uniformly in $\eta$, and this is sufficient 
to conclude that the set $\mathcal{I}$ of equilibrium measures of the environment process is a 
singleton, i.e., $\mathcal{I}=\{\mu_e\}$. Indeed, suppose that there are two equilibrium 
measures $\nu,\nu'\in\mathcal{I}$. Then
\be{uniqness}
\begin{aligned}
|\langle f \rangle_{\nu}-\langle f \rangle_{\nu'}| 
&=|\langle S(t)f \rangle_{\nu}-\langle S(t)f \rangle_{\nu'}|\\
&\leq |\langle S(t)f \rangle_{\nu}-\langle S(t)f \rangle_{\mu}|
+ |\langle S(t)f \rangle_{\nu'}-\langle S(t)f \rangle_{\mu}|\\
&= |\langle \big[S(t)f -\langle S(t)f \rangle_{\mu}]\rangle_{\nu}|
+ |\langle \big[S(t)f -\langle S(t)f\rangle_{\mu}]
\rangle_{\nu'}|\\
&\leq 2\left\| S(t)f -\langle S(t)f \rangle_\mu\right\|_\infty.
\end{aligned}
\ee
Since the l.h.s.\ of (\ref{uniqness}) does not depend on $t$, and the r.h.s.\ tends to zero as 
$t\to\infty$, we have $\nu=\nu'=\mu_e$. Next, $\mu_e$ is uniquely ergodic, meaning that the 
environment process converges to $\mu_e$ as $t\to\infty$ no matter what its starting distribution 
is. Indeed, for any $\mu'$,
\be{ergodicity}
|\langle S(t)f \rangle_{\mu'}-\langle S(t)f \rangle_{\mu}|
= | \langle \big[S(t)f -\langle S(t)f \rangle_{\mu}]\rangle_{\mu'}|
\leq \left\| S(t)f -\langle S(t)f \rangle_\mu\right\|_\infty,
\ee
and therefore 
\be{DOM}
\begin{aligned}
\langle f \rangle_{\mu_e}
&=\lim_{t\to\infty} S(t)f 
= \lim_{t\to\infty} \langle S(t)f \rangle_\mu 
= \lim_{t\to\infty} \left\langle \sum_{n\in\N} g_n(t,f) \right\rangle_\mu\\
&= \lim_{t\to\infty} \sum_{n\in\N}\langle g_n(t,f) \rangle_\mu
=  \sum_{n\in\N}\lim_{t\to\infty}\langle g_n(t,f) \rangle_\mu,
\end{aligned}
\ee
where the last equality is justified by the bound in (\ref{gSupBound}) in combination
with the dominated convergence theorem.
\end{proof}

We close this section by giving a more transparent description of $\mu_e$, more suitable
for explicit computation.

\bt{ErgodicLimit}
For $\alpha-\beta<\tfrac12 c$,
\be{FinalMeasure}
\langle f \rangle_{\mu_e} = \sum_{n\in\N} \langle\Psi_n\rangle_{\mu}
\ee
with
\be{RecursionPsi}
\Psi_1= f \quad \text{and}\quad \Psi_{n+1} 
= L_* L_0^{-1}(\Psi_n-\langle\Psi_n\rangle_{\mu}), \qquad n\in\N,
\ee
where $L_0^{-1} = \int_0^\infty S_0(t)\,\di t$ (whose domain is the set of all $f\in C(\Omega)$ 
with $\langle f\rangle_\mu=0$).
\et

\begin{proof}
By (\ref{DOM}), the claim is equivalent to showing that 
\be{Psi}
\lim_{t\to\infty} \langle g_n(t,f) \rangle_\mu = \langle\Psi_n\rangle_{\mu}.
\ee

First consider the case $n=2$. Then 
\be{Psi2}
\begin{aligned}
\lim_{t\to\infty} \langle g_2(t,f)\rangle_\mu
&=\lim_{t\to\infty} \left\langle \int_0^t \di s\,S_0(t-s)\,L_*\,g_1(s,f)\right\rangle_\mu\\
&= \lim_{t\to\infty} \left\langle \int_0^t \di s\,\,L_*\,g_1(s,f)\right\rangle_\mu\\
&= \lim_{t\to\infty} \left\langle \int_0^t \di s\,\,L_* S_0(s)f \right\rangle_\mu\\
&= \lim_{t\to\infty} \left \langle \int_0^t \di s\,\,L_* 
\big[S_0(s)(f-\langle f\rangle_{\mu})\big] \right\rangle_\mu\\
&= \left\langle\lim_{t\to\infty} L_*\int_0^t \di s\,\,
S_0(s)(f-\langle f\rangle_{\mu}) \right\rangle_\mu
= \langle L_* L_0^{-1}(f-\langle f\rangle_{\mu}) \rangle_\mu,
\end{aligned}
\ee
where the second equality uses that $\mu$ is invariant w.r.t.\ $S_0$, while the fifth equality
uses the linearity and continuity of $L_*$ in combination with the bound in (\ref{gSupBound}).

For general $n$, the argument runs as follows. First write
\be{Psin1}
\begin{aligned}
&\langle g_n(t,f)\rangle_\mu\\ 
&= \left\langle \int_0^t \di s\,S_0(t-t_1)\,L_*\,g_{n-1}(t_1,f) \right\rangle_\mu\\
&= \left\langle \int_0^t \di t_1\,\,L_*\,g_{n-1}(t_1,f) \right\rangle_\mu\\
&= \left\langle \int_0^t \di t_1 \int_0^{t_1} \di t_2 \cdots \int_0^{t_{n-1}}\di t_n\,
\big[L_* S_0(t_1-t_2) \cdots L_* S_0(t_{n-1}-t_n) L_* S_0(t_n)\big]\,f\right\rangle_\mu\\
&= \left\langle \int_0^t \di t_n \int_0^{t-t_n} \di t_{n-1} \cdots \int_0^{t-t_2}
\di t_1\,\big[L_* S_0(t_1) L_* S_0(t_2) \cdots L_* S_0(t_{n-1}) L_* S_0(t_n)\big]\,
f\right\rangle_\mu.
\end{aligned}
\ee
Next let $t\to\infty$ to obtain 
\be{Psin2}
\begin{aligned}
&\lim_{t\to\infty} \langle g_n(t,f)\rangle_\mu\\ 
&= \left\langle \int_0^\infty \di t_n\int_0^\infty\di t_{n-1} \cdots \int_0^\infty\di t_1\,
\big[L_* S_0(t_1) L_* S_0(t_2) \cdots L_* S_0(t_{n-1}) L_* S_0(t_n)\big]\,f\right\rangle_\mu\\
&= \left\langle L_*\int_0^\infty \di t_1\, S_0(t_1)\,L_*\int_0^\infty \di t_2\, S_0(t_2)
\cdots L_*\int_0^\infty \di t_n\, S_0(t_n)\,(f-\langle f\rangle_{\mu})\right\rangle_\mu\\
&=\left\langle L_*\int_0^\infty \di t_1\, S_0(t_1)\,L_*\int_0^\infty \di t_2\, S_0(t_2)
\cdots L_* L_0^{-1}(f-\langle f\rangle_{\mu})\right\rangle_\mu\\
&=\left\langle L_*\int_0^\infty \di t_1\, S_0(t_1)\,L_*\int_0^\infty \di t_2\, S_0(t_2)
\cdots L_*\int_0^\infty \di t_{n-1}\, S_0(t_{n-1})\Psi_2\right\rangle_\mu,
\end{aligned}
\ee
where we insert $L_*L_0^{-1}(f-\langle f\rangle_\mu)=\Psi_2$. Iteration shows that the latter 
expression is equal to
\be{Psin3}
\begin{aligned}
\left\langle L_*\int_0^\infty \di t_1\, S_0(t_1)\Psi_{n-1}\right\rangle_\mu
&= \left\langle L_*\int_0^\infty \di t_1\, 
S_0(t_1)(\Psi_{n-1}-\langle \Psi_{n-1}\rangle_{\mu})\right\rangle_\mu\\
&= \left\langle L_*L_0^{-1}(\Psi_{n-1}-\langle\Psi_{n-1}\rangle_{\mu})\right\rangle_{\mu} 
= \langle\Psi_n\rangle_{\mu}.
\end{aligned}
\ee
\end{proof}


\subsection{Expansion of the global speed}
\label{S3.3}

As we argued in (\ref{LLNlim}), the global speed of $X$ is given by 
\be{globsp}
v=(2\widetilde{\rho}-1)(\alpha-\beta)
\ee 
with $\widetilde{\rho}=\langle\phi_0\rangle_{\mu_e}$. By using Theorem~\ref{ErgodicLimit}, 
we can now expand $\widetilde{\rho}$.

First, if $\langle\phi_0\rangle_\mu = \rho$ is the particle density, then 
\be{Approx}
\widetilde{\rho} = \langle\phi_0\rangle_ {\mu_e} 
= \rho + \sum_{n=2}^\infty \langle \Psi_n\rangle_\mu,
\ee
where $\Psi_n$ is constructed recursively via (\ref{RecursionPsi}) with $f=\phi_0$. 
We have  
\be{PsiConstant} 
\langle\Psi_n\rangle_\mu = d_n\,(\alpha-\beta)^{n-1}, \quad n\in\N,
\ee
where $d_n=d_n(\alpha+\beta;P^\mu)$, and the factor $(\alpha-\beta)^{n-1}$ comes from 
the fact that the operator $L_*$ is applied $n-1$ times to compute $\Psi_n$, as is seen 
from (\ref{RecursionPsi}). Recall that, in (\ref{RevIrrev}), $L_\mathrm{SRW}$ caries the 
prefactor $\alpha+\beta$, while $L_*$ carries the prefactor $\alpha-\beta$. Combining 
(\ref{globsp}--\ref{Approx}), we have
\be{ApproxSpeed}
v = \sum_{n\in\N} c_n\,(\alpha-\beta)^n,
\ee
with $c_1=2\rho-1$ and $c_n=2d_n$, $n\in\N\backslash\{1\}$. 

For $n=2,3$ we have 
\be{2ndTerm}
\begin{aligned} 
c_2 &= 2 \big\langle\phi_0 L_0^{-1}\big(\phi_1-\phi_{-1}\big)\big\rangle_\mu\\
c_3 &= \tfrac12\left\langle\psi_0 L_0^{-1} \big[\psi_{-1} L_0^{-1}\bar\phi_{-2}
- \psi_{1}L_0^{-1}\bar\phi_0 - \psi_{-1}L_0^{-1}\bar\phi_0 + \psi_{1}L_0^{-1} 
\bar\phi_2\big]\right\rangle_\mu, 
\end{aligned}
\ee
where $\phi_i(\eta)=\eta(i)$, $\eta\in\Omega$, $\bar\phi_i=\phi_i-\langle\phi_i\rangle_{\mu}$ 
and $\psi_i=2\phi_i-1$. It is possible to compute $c_2$ and $c_3$ for appropriate choices 
of $\xi$.

If the law of $\xi$ is invariant under reflection w.r.t.\ the origin, then $\xi$ has the 
same distribution as $\xi'$ defined by $\xi'(x)=\xi(-x)$, $x\in\Z$. In that case $c_2=0$,
and consequently $v=(2\rho-1)(\alpha-\beta)+O((\alpha-\beta)^3)$. For examples of interacting 
particle systems with $M<\epsilon$, see Liggett~\cite{Li85}, Section I.4. Some of these 
examples have the reflection symmetry property.

An alternative formula for $c_2$ is (recall (\ref{RevIrrev})) 
\be{c2alt}
c_2 = 2 \int_0^\infty \di t\,
\Big(E_\mathrm{SRW,1}[K(Y_t,t)]-E_\mathrm{SRW,-1}[K(Y_t,t)]\Big),
\ee
where
\be{Kdef}
K(i,t) 
= E_{P^\mu}[\xi_0(0)\xi_t(i)]
= \langle\phi_0\,(S_\mathrm{IPS}(t)\phi_i)\rangle_\mu, \qquad i\in\Z,\,t\geq 0,
\ee 
is the space-time correlation function of the interacting particle system (with generator
$L_\mathrm{IPS}$), and $E_{\mathrm{SRW},i}$ is the expectation over simple random walk 
$Y=(Y_t)_{t\geq 0}$ jumping at rate $\alpha+\beta$ (with generator $L_\mathrm{SRW}$) 
starting from $i$. If $\mu$ is a \emph{reversible} equilibrium, then (recall (\ref{Pmuass}))
\be{revcalc}
K(i,t) = \langle\phi_0\,(S_\mathrm{IPS}(t)\phi_i)\rangle_\mu
= \langle(S_\mathrm{IPS}(t)\,\phi_0)\phi_i\rangle_\mu
= \langle(S_\mathrm{IPS}(t)\phi_{-i})\,\phi_0\rangle_\mu
= K(-i,t),
\ee
implying that $c_2=0$.

In Appendix~\ref{appB} we compute $c_3$ for the independent spin-flip dynamics, for which
$c_2=0$.


\appendix

\section{Examples of cone-mixing}
\label{appA}


\subsection{Spin-flip systems in the regime $M<\epsilon$}
\label{A.1}

Let $\xi$ be a spin-flip system for which $M<\epsilon$. We recall that in a 
spin-flip system only one coordinate changes in a single transition. The rate to 
flip the spin at site $x\in\Z$ in configuration $\eta\in\Omega$ is $c(x,\eta)$. 
As shown in Steif~\cite{St91} and in Maes and Shlosman~\cite{MaSc93}, two 
copies $\xi,\xi'$ of the spin-flip system starting from configurations $\eta,\eta'$
can be coupled such that, uniformly in $t$ and $\eta,\eta'$,
\be{Wasserstein}
\widehat{P}_{\eta,\eta'}\big(\exists\,s\geq t\colon\,\xi_s(x) \neq \xi'_s(x)\big)
\leq \sum_{ {y\in\Z:} \atop {\eta(y) \neq \eta'(y)} }  
\erom^{-\epsilon t}\,\big(\erom^{\Gamma t}\big)(y,x)
\leq \erom^{-(\epsilon-M)t}, 
\ee 
where $\widehat{P}_{\eta,\eta'}$ is the Vasershtein coupling (or basic coupling), and 
$\Gamma$ is the matrix $\Gamma =(\gamma(u,v))_{u,v\in\Z}$ with elements
\be{Gamma}
\gamma(u,v) = \sup_{\eta\in\Omega} |c(u,\eta)-c(u,\eta^v)|. 
\ee
Recall (\ref{Mepsdef}) to see that $\Gamma$ is a bounded operator on $\ell_1(\Z)$ with 
norm $M$ (see also Liggett~\cite{Li85}, Section I.3).

Define
\be{rhotdef}
\rho(t) = \sup_{\eta,\eta'\in\Omega} \widehat{P}_{\eta,\eta'}\big(\exists\,s \geq t\colon\,
\xi_s(0) \neq \xi'_s(0)\big), \qquad t\geq 0.
\ee
Recall Definition~\ref{PhiMixing}, fix $\theta\in(0,\tfrac12\pi)$ and put $c=c(\theta)= 
\cot\theta$. For $B\in\cF_t^\theta$, estimate
\be{ConeM<e}
\begin{aligned}
|P_{\eta}(B)-P_{\eta'}(B)| 
&\leq \widehat{P}_{\eta,\eta'}
\big(\exists\,x\in\Z\,\,\exists\,s \geq t+c|x|\colon\,\xi_s(x)\neq \xi'_s(x)\big)\\
&\leq \sum_{x\in\Z} \widehat{P}_{\eta,\eta'}
\big(\exists\,s \geq t\colon\,\xi_s(x) \neq \xi'_s(x)\big)\\
&\leq \sum_{x\in\Z} \rho(t+c|x|)\\
& \leq \rho(t) + 2 \int_0^\infty \rho(t+cu)\,\di u\\
&= \rho(t) +\frac{2}{c} \int_0^\infty \rho(t+v)\,\di v.  
\end{aligned}
\ee
Since this estimate is uniform in $B$ and $\eta,\eta'$, it follows that for the cone mixing
property to hold it suffices that
\be{conesuf}
\int_0^\infty \rho(v)\,\di v <\infty.
\ee
It follows from (\ref{Wasserstein}) that $\rho(t) \leq \erom^{-(\epsilon-M)t}$, which indeed 
is integrable. 

Note that if the supremum in (\ref{rhotdef}) is attained at the same pair of starting 
configurations $\eta,\eta'$ for all $t \geq 0$, then (\ref{conesuf}) amounts to the condition 
that the average coupling time at the origin for this pair is finite.


\subsection{Attractive spin-flip dynamics}
\label{A.2}

An attractive spin-flip system $\xi$ has rates $c(x,\eta)$ satisfying
\be{attractivness}
\begin{aligned}
&c(x,\eta) \leq c(x,\eta') \quad \text{ if } \eta(x)=\eta'(x)=0,\\
&c(x,\eta) \geq c(x,\eta') \quad \text{ if } \eta(x)=\eta'(x)=1,
\end{aligned}
\ee 
whenever $\eta\leq\eta'$ (see Liggett~\cite{Li85}, Chapter III). If $c(x,\eta) = c(x+y,\tau_y\eta)$
for all $y\in\Z$, then attractivity implies that, for any pair of configurations $\eta,\eta'$,
\be{+-coupling}
\widehat{P}_{\eta,\eta'}\big(\exists\,s \geq t\colon\,\xi_s(x) \neq \xi'_s(x)\big)
\leq \widehat{P}_{[0],[1]}\big(\exists\,s \geq t\colon\,\xi_s(0)\neq \xi'_s(0)\big),
\ee 
where $[0]$ and $[1]$ are the configurations with all $0$'s and all $1$'s, respectively. 
Proceeding as in (\ref{ConeM<e}), we find that for the cone-mixing property to hold it 
suffices that
\be{conesufalt}
\int_0^\infty \rho^*(v)\,\di v < \infty, \qquad
\rho^*(t) = \widehat{P}_{[0],[1]}\big(\exists\,s \geq t\colon\,\xi_s(0)\neq \xi'_s(0)\big). 
\ee 

Examples of attractive spin-flip systems are the (ferromagnetic) Stochastic Ising Model, the 
Contact Process, the Voter Model, and the Majority Vote Process (see Liggett~\cite{Li85}, 
Chapter III). For the one-dimensional Stochastic Ising Model, $t\mapsto\rho^*(t)$ decays 
exponentially fast at any temperature (see Holley~\cite{Ho85}). The same is true for the
one-dimensional Majority Vote Process (Liggett~\cite{Li85}, Example~III.2.12). Hence both
are cone-mixing. The one-dimensional Voter Model has equilibria $p\delta_{[0]}+(1-p)\delta_{[1]}$, 
$p\in [0,1]$, and therefore is not interesting for us. The Contact Process has equilibria
$p\delta_{[0]}+(1-p)\nu$, $p \in [0,1]$, but $\nu$ is not cone-mixing. 

In view of the remark made at the end of Section \ref{S1.4}, we note the following. For the 
Stochastic Ising Model in dimensions $d \geq 2$ exponentially fast decay occurs only at high 
enough temperature (Martinelli~\cite{Ma97}, Theorem 4.1). The Voter Model in dimensions 
$d \geq 3$ has non-trivial ergodic equilibria, but none of these is cone-mixing. The same 
is true for the Contact Process in dimensions $d \geq 2$.


\subsection{Space-time Gibbs measures}
\label{A.3} 

We next give an example of a discrete-time dynamic random environment that is cone-mixing 
but not Markovian. Accordingly, in (\ref{Fdefs}) we must replace $\cF_0$ by $\cF_{-\N_0} 
= \{\xi_t(x)\colon\,x\in\Z,\,t\in (-\N_0)\}$. Let $\sigma=\{\sigma(x,y)\colon\,(x,y)\in\Z^2\}$ 
be a two-dimensional Gibbsian random field in the Dobrushin regime (see Georgii \cite{Ge88}, 
Section 8.2). We can define a discrete-time dynamic random environment $\xi$ on $\Omega$ 
by putting
\be{2Gibbs}
\xi_t(x) = \sigma(x,t) \qquad (x,t)\in\Z^2.
\ee 
The cone-mixing condition for $\xi$ follows from the mixing condition of $\sigma$ in the 
Dobrushin regime. In particular, the decay of the mixing function $\Phi$ in (\ref{conecor}) 
is like the decay of the Dobrushin matrix, which can be polynomial.


\section{Independent spin-flips}
\label{appB}

Let $\xi$ be the Markov process with generator $L_\mathrm{ISF}$ given by
\be{GenISF}
(L_\mathrm{ISF}f)(\eta) = \sum_{x\in\Z} c(x,\eta)\,\big[f(\eta^x)-f(\eta)\big],
\qquad \eta\in\Omega,
\ee 
where 
\be{FlipRate}
c(x,\eta) = \gamma[1-\eta(x)] + \delta\eta(x),
\ee
i.e., 0's flip to 1's at rate $\gamma$ and 1's flip to 0's at rate $\delta$, independently
of each other. Such a $\xi$ is an example of a dynamics with $M<\epsilon$, for which 
Theorem~\ref{ErgodicLimit} holds. From the expansion of the global speed in (\ref{ApproxSpeed}) 
we see that $c_2=0$, because the dynamics is invariant under reflection in the origin. We 
explain the main ingredients that are needed to compute $c_3$ in (\ref{third}).

The equilibrium measure of $\xi$ is the Bernoulli product measure $\nu_\rho$ with parameter 
$\rho=\gamma/(\gamma+\delta)$. We therefore see from (\ref{2ndTerm}) that we must compute 
expressions of the form
\be{1c_3}
I(j,i) = \left\langle (2\eta(0)-1) L_0^{-1} 
\big[(2\eta(j)-1)L_0^{-1}(\eta(i)-\rho)\big]\right\rangle_{\nu_\rho},
\ee
where $\eta$ is a typical configuration of the environment process $\zeta=(\zeta_t)_{t\geq 0} 
= (\tau_{X_t}\xi_t)_{t\geq 0}$ (recall Definiton~\ref{EnvPRo}), and  
\be{A}
(j,i)\in A = \{(-1,-2),(-1,0),(1,0),(1,2)\}.
\ee
By Lemma~\ref{GenDecomposition} we have $L_0=L_\mathrm{SRW}+L_\mathrm{ISF}$, with $L_\mathrm{SRW}$
the generator of simple random walk on $\Z$ jumping at rate $U=\alpha+\beta$. Hence
\be{2c_3}
(S_0(t)\eta)(i) = E_R^\eta[\eta_t(i)]
= \sum_{y\in\Z} p_{Ut}(0,y)\,E_\mathrm{ISF}^{\tau_y\eta}[\eta_t(i)]
= \sum_{y\in\Z} p_{Ut}(0,y)\,E_\mathrm{ISF}^\eta[\eta_t(i-y)],
\ee
where $\tau_y$ is the shift of space over $x$,
\be{3c_3}
E_\mathrm{ISF}^\eta[\eta_t(i)] =
\eta(i)\,\erom^{-Vt} + \rho(1-\erom^{-Vt})
\ee
with $V=\gamma+\delta$, and $p_t(0,y)$ is the transition kernel of simple random walk on $\Z$ 
jumping at rate $1$. Therefore, by (\ref{2c_3}--\ref{3c_3}),
we have 
\be{4c_3}
L_0^{-1}(\eta(i)-\rho)
= \int_0^\infty S_0(t)(\eta(i)-\rho)\,\di t
= \sum_{y\in\Z} \eta(i-y)\,G_V(y)-\rho\,\frac{1}{V}
\ee
with 
\be{G_c}
G_V(y) = \int_0^\infty \erom^{-Vt}\,p_{Ut}(0,y)\,\di t.
\ee

With these ingredients we can compute (\ref{1c_3}), ending up with
\be{5c_3}
c_3 = \sum_{(j,i)\in A} I(j,i)
= \frac{4}{U}\,\rho(2\rho-1)(1-\rho)\,
\left[\frac{2U+V}{U}\,G_V(0) - \frac{3U+2V}{U}\,G_{2V}(0) - G_{2V}(1)\right].
\ee
The expression between square brackets can be worked out, because 
\be{G_c(0)} 
G_V(0) = \int_0^\infty \erom^{-Vt}\,p_{Ut}(0,0)\,\di t
= \frac{1}{2\pi} \int_{-\pi}^{\pi} 
\frac{\di \theta}{(U+V)-U \cos\theta}
= \frac{1}{\sqrt{(U+V)^2-U^2}}
\ee
and 
\be{G1-0}
G_V(1) = \frac{U+V}{U}\,G_V(0)- \frac{1}{U},
\ee
where the latter is derived by using that
\be{6c_3} 
\frac{\partial}{\partial t}\,p_{Ut}(0,0)
=\tfrac12 U\,\big[p_{Ut}(0,1)+p_{Ut}(0,-1)-2p_{Ut}(0,0)\big]
\ee  
and $p_{Ut}(0,1)=p_{Ut}(0,-1)$. This leads to (\ref{third}).


\end{document}